\documentclass[review,hidelinks,onefignum,onetabnum]{siamart250211}

\nolinenumbers

\usepackage{lipsum}
\usepackage{amsfonts}
\usepackage{graphicx}
\usepackage{epstopdf}
\usepackage{algpseudocode}
\ifpdf
  \DeclareGraphicsExtensions{.eps,.pdf,.png,.jpg}
\else
  \DeclareGraphicsExtensions{.eps}
\fi

\usepackage[normalem]{ulem}
\usepackage{subfigure}
\usepackage{amssymb,amsfonts,amsmath,latexsym,dsfont}
\usepackage{mathtools}
\usepackage{wrapfig}

\usepackage{tikz,pgfplots}
\usetikzlibrary{plotmarks}

\usepackage{graphicx}
\usepackage{mathrsfs}
\usepackage{framed,multirow}
\usepackage{setspace}

\usepackage{booktabs}

\usepackage{soul}

\newdimen\iwidth
\newdimen\iheight

\newcommand{\mc}[3]{\multicolumn{#1}{#2}{#3}}

\newtheorem*{remark}{Remark}

\graphicspath{{Fig/}{./}}

\definecolor{rachelis}{rgb}{0.08, 0.25, 0.68}

\definecolor{yunhuis}{rgb}{0.58, 0.08, 0.35}

\definecolor{erans}{rgb}{0.78, 0.15, 0.08}

\usepackage[textsize=tiny]{todonotes}

\newcommand{\im}{\textit{\i}}

\newcommand{\bfe}{{\bf e}}

\newcommand{\bfx}{ {\bf x}}

\newcommand{\bfu}{{\bf u}}
\newcommand{\bfq}{{\bf q}}

\newcommand{\bfr}{{\bf r}}

\newsiamremark{hypothesis}{Hypothesis}
\crefname{hypothesis}{Hypothesis}{Hypotheses}
\newsiamthm{claim}{Claim}
\newsiamremark{fact}{Fact}
\crefname{fact}{Fact}{Facts}

\headers{Vanka-Smoothed Multigrid for Helmholtz}{R Yovel, Y. He, and E. Treister}

\title{Vanka-smoothed Shifted Laplacian multigrid preconditioners for the Helmholtz equations\thanks{Submitted to the editors November 18th 2025.
\funding{This research was supported by The Israel Science Foundation (grant No. 656/23). RY is supported by the Ariane de Rothschild scholarship and by Kreitman Negev-Faran scholarship. The authors also thank the Lynn and William Frankel Center at the Stein faculty for Computer and Information Science at BGU.}}}

\author{Rachel Yovel\thanks{Institute for Interdisciplinary Computational Science, the Stein Faculty of Computer and Information Science, Ben-Gurion University of the Negev, Beer-Sheva, Israel
  (\email{yovelr@bgu.ac.il}, \email{erant@bgu.ac.il}).}
\and Yunhui He\thanks{Department of Mathematics, Houston University, Houston, Texas, USA
  (\email{yhe43@central.uh.edu}}).
\and Eran Treister$^\dag$}

\usepackage{amsopn}

\ifpdf
\hypersetup{
  pdftitle={Vanka-Smoothed Multigrid for Helmholtz},
  pdfauthor={R. Yovel, Y. He, and E. Treister}
}
\fi

\begin{document}

\maketitle

\begin{abstract}
We present an improved multigrid preconditioner for the acoustic Helmholtz equation with enhanced scalability. 
Standard multigrid fails to converge for the Helmholtz equation, and the well-known complex shifted Laplacian method overcomes it by adding a complex shift and using the shifted system as a preconditioner.
However, the added complex shift grows with the frequency and interferes with the preconditioner's scalability.
In this work, we present an additive Vanka smoother that requires a much lower shift than point-wise smoothers, and thereby enhances the scalability.
By carefully designing different ingredients of the multigrid cycle, the presented method enables deep V-cycles with a small and bounded shift, even when many levels are used. 
We validate our method theoretically by local Fourier analysis, and hold numerical experiments for homogeneous and heterogeneous media. 
We show that our method outperforms plain shifted Laplacian in terms of runtimes and performs well on heterogeneous media in 2D and challenging geophysical media in 3D.
\end{abstract}



\section{Introduction}

The acoustic Helmholtz equation models wave propagation in the frequency domain, and is given by 
\begin{equation} \label{eq:acousitcHelm}
-\Delta p - \omega^{2}\kappa^2\left(1-\frac{\gamma}{\omega}\im\right)p = q
\end{equation}
where $p = p(\vec{x}), \, \vec{x}\in\Omega$ is the Fourier transform of the wave's pressure field, $\omega = 2\pi f$ is the angular frequency, where $f$ denotes the temporal frequency of the wave, $\kappa = \kappa(\vec{x}) > 0$ is the ``slowness'' of the wave in the medium (the inverse of the wave velocity), and $q(\vec{x})$ is the source of the waves. The notation $\im$ stands for the imaginary unit and $\gamma$ represents the physical attenuation.


We discretize the acoustic Helmholtz equation \eqref{eq:acousitcHelm} by a finite-difference scheme on a regular grid, typically equipped with absorbing boundary conditions (ABC) \cite{engquist1977absorbing} or perfectly matched layers (PML)
\cite{berenger1994perfectly,   harari2000analytical, rabinovich2010comparison}, to mimic the propagation of a wave in an open domain, and to avoid reflections from the boundary. 
The resulting linear system is large, indefinite and complex-valued due to the ABC and possible attenuation. 
High-frequency Helmholtz problems require very fine meshes and hence result in a large number of unknowns \cite{bayliss1985accuracy,haber2011fast}. Solving the discretized equation at large scale 3D scenarios is challenging, and is still considered as open problem.


Many preconditioners have been suggested for the Helmholtz equation, most of which are based on domain decomposition methods \cite{gander2013domain, stolk2013rapidly, tsuji2014sweeping, chen2016robust, dolean2025achieving}
and multigrid methods \cite{erlangga2006novel, olson2010smoothed, oosterlee2010shifted, livshits2014scalable, tsuji2015augmented, cools2015multi}.
The well-known complex-shifted Laplacian multigrid preconditioner (CSLP) \cite{erlangga2006novel} suggests adding a frequency-dependent complex shift to the Helmholtz equation \eqref{eq:acousitcHelm}, and using the shifted version (solved by multigrid) as a preconditioner for the original system. 
Adding a sufficiently large complex shift promises multigrid convergence \cite{elman2001multigrid}. 
Because the mesh must be refined as the frequency increases, the corresponding growth of the added shift reduces the preconditioner’s efficiency for large wavenumbers, making the method non-scalable with respect to grid size.
This gives rise to seeking multigrid methods for the Helmholtz equation that require less shift.


Over the years, many works have inquired as to what extent the shift can be reduced while maintaining good convergence properties.
For a two level deflation method, a fully shiftless and scalable method was presented in \cite{Dwarka2020}, using bicubic intergrid operators \cite{donatelli2011grid}, rather than the standard bilinear interpolation and full weighting restriction \cite{trottenberg2000multigrid}. Nevertheless, the coarse problem in a two-level method might be too large and this scalability is harder to achieve for a multilevel method. 
In \cite{chen2024matrix}, the two-level deflation method of \cite{Dwarka2020} was adapted for parallel matrix-free implementation, and equipped with a multigrid method as a coarse solver.
A multigrid method with no shift on the fine level was presented in \cite{cools2014new}.
The method uses a level-dependent shift strategy, leaving the fine level unmodified while adjusting the coarser levels.
However, the increasing added shift in the coarser levels prevents wavenumber independent convergence.
In \cite{chen2025matrix}, a scalable multilevel deflation method was presented for Helmholtz problems. 
However, the reported scalability for more than two levels was obtained there for moderate to low wavenumbers (e.g., corresponding to approximately 20 grid points per wavelength) and a limited number of levels.
Achieving scalability in more challenging regimes with fewer grid points per wavelength remains an open issue.


When the number of levels is limited, the computational bottleneck of multigrid methods lies in the coarse grid problem, which may still be too large to solve directly.
Many special adaptations were suggested to deal with the coarse grid problem in large 3D instances:
In \cite{treister2024hybrid}, domain decomposition was used as a coarse solver, and in \cite{yovel2024lfa}, the coarse grid was solved iteratively using a Kaczmarz relaxation \cite{gordon2013robust,li20152d}.
In \cite{oosterlee2000krylov}, in the context of fluid flow, a multigrid method is accelerated by using an iterative approach on a coarser grid.
More generally, in extremely large problems the coarse grid solve becomes a dominant cost and a key limiting factor, due to loss of sparsity and the superlinear complexity growth of LU factorization.


The pair of works \cite{gander2015applying, cocquet2017large} addresses the two complementary questions: how small the shift must be to retain the scalability of CSLP, and how large it should be to promise multigrid convergence. 
On the one hand, it was shown in \cite{gander2015applying} that in order to achieve scalability, the shift must grow at most linearly with the wavenumber.
On the other hand, it was shown in \cite{cocquet2017large} that in order to ensure multigrid convergence the shift must grow super-linearly with the wavenumber (and even, in the general case, quadratically). Consequently, there is very little hope for the development of a fully scalable multigrid method for the Helmholtz equation.

Our aim in this work is to reduce the shift as much as possible, and to improve scalability within these limitations.
To address the first problem we suggest a careful design of intergrid operators that lowers the shift size while keeping the operator complexity of the Galerkin coarse approximation bounded.
As high order intergrid results in a wide computational stencil for the Galerkin operator, we suggest a level-dependent intergrid scheme that prioritizes the first level.
That is, we take high-order intergrid between the first and the second level, and mixed-order intergrid between the rest of the levels.
To enable cycles with a large number of levels without a significant increase in computational cost, we carefully design an additive Vanka relaxation,
and show that it achieves  better scalability, using a relatively small Vanka patch thus keeping a reasonable computational effort per iteration.
For a given grid size and wavenumber in homogeneous or mildly heterogeneous media, the shift, iteration count and time do not grow using a certain number of levels and above.
In challenging geophysical media the time grows mildly with the level, and a suggested combination of damped Jacobi and Vanka smoothers enables to control this growth.
While the shift in this method still grows quadratically with the wavenumber, the iteration count is significantly lower, since the shifts are smaller.
Moreover, the designed multigrid components enable deep V-cycles, which are substantially cheaper compared to deep W-cycles.
The multiple coarse-grid visits in W-cycles make it possible to solve the Helmholtz equation even with a standard (and potentially unstable) smoother.
However, this becomes impractical when the number of levels increases. Our method provides a remedy by convergent V-cycles equipped with block-smoothers.

We support our method with local Fourier analysis (LFA) and validate it by numerical experiments for various homogeneous and heterogeneous Helmholtz problems in 2D and 3D.
Overall, the combination of the level-dependent intergrid scheme with the tailored Vanka patch facilitates the solution of large 3D problems in real-world geophysical media, yielding  relatively low iteration counts and improved time scalability.

The remainder of the paper is organized as follows: in Section \ref{sec:background} we give background on multigrid methods, with emphasis on complex-shifted Laplacian methods for Helmholtz equations and LFA for additive Vanka smoothers.
In Section \ref{sec:method} we present our method, in Section \ref{sec:lfa} we support our method by LFA and in Section \ref{sec:results} we give numerical results.
Finally, in Section \ref{sec:conclusion} we make concluding remarks.

\section{Background} \label{sec:background}

\subsection{Discretization}\label{subsec:discretization}

The acoustic Helmholtz equation \eqref{eq:acousitcHelm} is typically solved on a bounded domain, discretized by a regular grid with a finite differences stencil.
However, the standard second order central difference scheme (which yields a 5-point stencil in 2D and a 7-point stencil in 3D) might lead to large numerical dispersion, and hence to poor multigrid convergence \cite{stolk2014multigrid, yovel2024lfa}. 

A compact 4-th order stencil for the Helmholtz equation was suggested in  \cite{singer1998high}
\begin{equation} \label{eq:disc4thCompact}
H = \frac{1}{h^2}
\begin{bmatrix}
-\frac{1}{6} & -\frac{2}{3}  & -\frac{1}{6} \\[0.4em]
-\frac{2}{3}  & \frac{10}{3}  & -\frac{2}{3}  \\[0.4em]
-\frac{1}{6} & -\frac{2}{3}  & -\frac{1}{6} 
\end{bmatrix}
- \kappa^2 \omega^2 \left(1-\frac{\gamma}{\omega} \im\right)
\begin{bmatrix}
 &  \frac{1}{12} & \\[0.4em]
 \frac{1}{12} & \frac{2}{3} &  \frac{1}{12} \\[0.4em]
 & \frac{1}{12} & 
\end{bmatrix},
\end{equation}
and later validated in the context of multigrid in \cite{umetani2009multigrid}.
In 3D, the corresponding compact $3\times3\times3$ stencil reads
\begin{equation}\label{eq:disc4thCompact3D}
H = L - \kappa^2 \omega^2 \left(1-\frac{\gamma}{\omega} \im\right) M,
\end{equation}
where
\begin{equation}
L = -\frac{1}{6h^2}  
\begin{bmatrix}
\begin{bmatrix}
 & 1 & \\
 1 & 2 & 1 \\
 &  1 &
\end{bmatrix}
\begin{bmatrix}
 1 & 2 & 1  \\
2 & -24 & 2  \\
1 & 2 & 1
\end{bmatrix}
\begin{bmatrix}
 & 1 & \\
 1 & 2 & 1 \\
 &  1 &
\end{bmatrix}
\end{bmatrix}
\end{equation}
and
\begin{equation}
M = \frac{1}{12}
\begin{bmatrix}
\begin{bmatrix}
 &  & \\
 & 1 &  \\
 &   &
\end{bmatrix}
\begin{bmatrix}
 & 1 &   \\
1 & 6 & 1  \\
  & 1 &  
\end{bmatrix}
\begin{bmatrix}
 &  & \\
 & 1 &  \\
 &   &
\end{bmatrix}
\end{bmatrix}.
\end{equation}

Throughout this work, we discretize \eqref{eq:acousitcHelm} on a nodal grid using the stencil \eqref{eq:disc4thCompact} in 2D or \eqref{eq:disc4thCompact3D} in 3D.

\subsection{Multigrid methods} \label{subsec:multigrid}

Multigrid \cite{brandt1977multi, trottenberg2000multigrid} is a family of iterative methods for the solution of linear systems that arise, typically, from the discretization of elliptic partial differential equations (PDEs).
The idea behind multigrid methods is based on the observation that standard iterative methods, such as damped Jacobi or Gauss-Seidel, tend to reduce the oscillatory (high-frequency) error components. 
As a complementary process, the \emph{coarse-grid correction} --- the use of the exact error on a coarser grid to estimate and correct the error on the fine grid --- reduces the smooth (low-frequency) error modes.

For a more detailed description of {a multigrid algorithm, let 
\begin{equation}\label{eq:linsys}
A_h\bfu=\bfq
\end{equation}
be the system we aim to solve, where $A_h$ is a discretized version of a given operator on a fine grid. 
Let $P$ be the prolongation operator that interpolates coarse grid vectors to the fine grid, and let $R$ be the restriction that gives a weighted sampling of a fine grid vector on the coarse grid. 
Let $A_H$ be the \emph{coarse grid operator} that approximates $A_h$ on the coarse grid. 
The coarse-grid correction can be calculated by solving the equation
\begin{equation*}
A_H\bfe_H = \bfr_H = R(\bfq-A_h\bfu^{(k)})
\end{equation*}
where $\bfe_H$ and $\bfr_H$ are the error and the residual on the coarse grid, respectively, and interpolating the solution back to the fine grid afterward:
\begin{equation*}
\bfe = P\bfe_H.
\end{equation*}
The two-grid cycle is summarized in
Algorithm \ref{alg:TwoCycle}. 
In matrix form, the two-grid error propagation operator is given by:
\begin{equation}\label{eq:2G}
TG = S^{\nu_2}(I-P A_H^{-1} R A_h)S^{\nu_1}
\end{equation}
where $S$ is the smoother's error propagation matrix, and $\nu_1,\nu_2$ are the number of pre- and post relaxations, respectively.

\begin{algorithm}
\caption{Two-grid cycle}\label{alg:TwoCycle}
\begin{algorithmic}[1]
\Statex \textbf{Algorithm:} $\bfu\leftarrow TwoGrid(A_h,\bfq,\bfu).$
\State Apply $\nu_1$ pre-relaxations: $\mathbf{u} \gets \text{Relax}(A_h,\mathbf{u},\mathbf{q})$
\State Compute and restrict the residual: $\mathbf{r}_H \gets R(\mathbf{q} - A_h\mathbf{u})$
\State Compute $\mathbf{e}_H$ by solving the coarse-grid problem $A_H\mathbf{e}_H = \mathbf{r}_H$
\State Apply coarse-grid correction: $\mathbf{u} \gets \mathbf{u} + P\mathbf{e}_H$
\State Apply post-relaxations: $\mathbf{u} \gets \text{Relax}(A_h,\mathbf{u},\mathbf{q})$
\end{algorithmic}
\end{algorithm}

The multigrid cycle is comprised of recursive application of Algorithm~\ref{alg:TwoCycle}. 
To obtain a multigrid V-cycle, step 3 in Algorithm \ref{alg:TwoCycle} is replaced by the entire algorithm on a coarser level, and this is repeated recursively until reaching the desired number of levels.
Similarly, to obtain a W-cycle, step 3 in Algorithm \ref{alg:TwoCycle} is replaced by two applications of the entire algorithm on a coarser grid.
It is worth mentioning that in W-cycles, the number of coarse grid solutions per cycle grows exponentially with the number of levels, while in V-cycles, the coarse grid problem is solved only once in a cycle.
It implies that deep W-cycles are substantially more expensive compared to deep V-cycles, if the solution of the coarsest grid is not trivial.
For Helmholtz problems, as we use more levels, we increase the relative number of grid points per wavelength for the coarsest grid, and hence standard cycles tend to diverge without a special care.
Our main goal in this work is to develop a deep V-cycle that manages to solve the Helmholtz equation efficiently.

\subsection{Shifted Laplacian multigrid}

Standard multigrid methods do not converge for Helmholtz problems, since standard smoothers are not stable for indefinite systems, and the coarse-grid correction can increase the error \cite{elman2001multigrid}.
The complex-shifted Laplacian preconditioner (CSLP) \cite{erlangga2006novel} 
suggests adding a complex shift to the acoustic Helmholtz equation \eqref{eq:acousitcHelm}, such that the shifted version can be solved using multigrid, and serves as a preconditioner for the original version.
Physically speaking, additional complex shift in the frequency domain is equivalent to additional attenuation in the time space, so that the waves decay faster.

The shifted operator is defined by
\begin{equation}\label{eq:shift}
H_s = H - \im\alpha\omega^2M,
\end{equation}
where $\alpha$ is a shifting parameter.
Given a sufficient shift $\alpha$, standard multigrid converges for $H_s$, and enables using it as a preconditioner for \eqref{eq:acousitcHelm} inside a Krylov method such as (flexible) GMRES \cite{saad1993flexible} or BiCGSTAB \cite{van1992bi}.

As shown in \cite{cocquet2017large}, for a convergent multigrid method, the added complex shift in \eqref{eq:shift} should be proportional to $\omega^2$. 
However, a larger shift results in $H_s$ that is not a good approximation of $H$, and hence, shifted Laplacian methods are not scalable with respect to the wavelength.
In fact, it was shown in \cite{gander2015applying} that if the shift grows super-linearly in $\omega$, a wavenumber independent convergence cannot be achieved. 
Hence, lowering the value of $\alpha$, or reducing the dependence of $\alpha$ on the depth of the multigrid cycle, is the best that can be done to improve the scalability properties of current shifted Laplacian methods.

\subsection{Local Fourier analysis}

LFA \cite{trottenberg2000multigrid} is an analytical tool used to predict the convergence rate of multigrid cycles \cite{brandt1977multi}. Under the assumption of periodic BC, Fourier modes are eigenfunctions of discretized differential operators.
LFA uses Fourier symbol calculation of an operator with periodic BC to approximately predict the convergence rate of multigrid methods for a corresponding operator with any other BC.
Smoothing LFA determines the Fourier symbol of the smoother, and is used to predict the convergence rate of the entire cycle.
For indefinite problems such as the Helmholtz equation, the coarse-grid correction is not ideal \cite{trottenberg2000multigrid} and smoothing LFA gives an overly optimistic prediction.
In this case, two-grid LFA should be considered, as it takes into account both the smoother and the coarse-grid correction.

\begin{definition}[\cite{trottenberg2000multigrid}, Chapter 8] \label{def_mu_loc}
    Let $\widetilde{S}(\theta)$ be the symbol matrix of the error propagation matrix $S$, for $\theta\in\left[-\frac{\pi}{2},\frac{3\pi}{2}\right]^2$, and let $T^\mathrm{high}=\left[-\frac{\pi}{2},\frac{3\pi}{2}\right]^2 \setminus \left[-\frac{\pi}{2},\frac{\pi}{2}\right]^2$.
Then the LFA smoothing factor is
    \[
        \mu_\mathrm{loc} \coloneqq \sup_{\theta\in T^\mathrm{high}} \rho(\widetilde{S}(\theta)).
    \]
\end{definition}

For overlapping smoothers, and particularly for the Vanka smoother, the calculation of $\widetilde{S}$ is not straightforward, and special techniques are required to perform smoothing analysis. 
Such techniques for multiplicative smoothers are given in \cite{sivaloganathan1991use, maclachlan2011local, rodrigo2016local, treister2024hybrid}, and for additive Vanka smoothers in \cite{farrell2021local, greif2023closed}.

\begin{definition}[\cite{trottenberg2000multigrid}, Chapter 8] \label{def_rho_loc}
    Let $\widetilde{TG}(\theta)$ be the symbol of the two-grid error operator \eqref{eq:2G} and let $T^\mathrm{low}=\left[-\frac{\pi}{2},\frac{\pi}{2}\right]^2$. Then the two-grid LFA convergence factor is defined by
\begin{equation}\label{eq:rho_loc}
\rho_\mathrm{loc} \coloneqq \sup_{\theta\in T^\mathrm{low}} \rho(\widetilde{TG}(\theta)).
\end{equation}
\end{definition}

To calculate the Fourier symbol of the two-grid error operator \eqref{eq:2G}, a space of harmonics should be introduced because Fourier modes are not eigenfunctions of the restriction $R$ and the prolongation $P$, not even when assuming periodic BC, because of the aliasing phenomenon. In 2D, e.g., four different frequencies alias to the same low frequency $\theta\in T^\mathrm{low}$, which gives rise to the definition of the 4-dimensional space of harmonics:
\begin{definition}[\cite{trottenberg2000multigrid}, Chapter 8] \label{def:4harmonics}
The 4-dimensional space of harmonics for $\theta$ is
\begin{equation}
E(\theta) = \text{span} \left\{ \theta, \theta', \theta'', \theta''' \right\}
\end{equation}
where $\theta' = [\theta_1+\pi, \theta_2+\pi]^T$, $\theta'' = [\theta_1+\pi, \theta_2]^T$, and $\theta''' = [\theta_1, \theta_2+\pi]^T$.
\end{definition}

The symbol matrix of the two-grid error operator is hence calculated by:
\begin{equation} \label{eq:symbolTG}
\widetilde{TG}(\theta) = \widetilde{\mathbf{S}}(\theta)^{\nu_2} (I-\widetilde{P}(\theta) \widetilde{H}_c^{-1}(\theta) \widetilde{R}(\theta) \widetilde{H}(\theta)) \widetilde{\mathbf{S}}(\theta)^{\nu_1}.
\end{equation}
A detailed two-grid LFA for the acoustic Helmholtz equation is described in \cite{cools2013local}.
As elaborated there, for the 2D case, $\widetilde{\mathbf{S}}(\theta)$ and $\widetilde{H}(\theta)$ are $4\times4$ diagonal symbol matrices, with the corresponding symbol evaluated in each of the four harmonics, and $\widetilde{R}(\theta)$ and $\widetilde{P}(\theta)$ are $1\times 4$ and $4\times 1$ matrices, similarly.

\section{Proposed method} \label{sec:method}

Our method is based on carefully designing the different multigrid components, to lower the required shift in the shifted Laplacian multigrid. 
In Subsection \ref{subsec:intergrid} we describe the coarse-grid correction, comprised of a level-dependent intergrid scheme and Galerkin coarse-grid correction.
In subsection \ref{subsec:smoother} we present the additive Vanka smoother and present different patches in 2D and 3D.

\subsection{Coarse-grid correction} \label{subsec:intergrid}

To design an efficient multigrid method, the coarse grid operator should resemble the fine grid operator, without overly increasing the operator complexity.
Generally speaking, the coarse grid operator $H_c$ can be calculated either by the Galerkin coarse approximation $H_c=R H P$, or by re-discretizing the PDE on a coarser grid. 
In this paper we use Galerkin coarse approximation, since it has theoretical advantages as a projection operator, and shows good convergence in practice for Helmholtz problems \cite{treister2024hybrid}.

Since the Galerkin operator is defined by the intergrid operators, its efficiency and operator complexity depend on the choice of intergrid operators. 
On the one hand,
high-order intergrid operators yield a better coarse-grid correction \cite{cools2014new, Dwarka2020}, but can increase the operator complexity.
On the other hand, standard intergrid operators (such as bilinear interpolation and full-weighting restriction) keep the operator complexity low but provide a weaker coarse-grid approximation of the fine-grid Helmholtz operator as the number of levels increases.
The balance between convergence and complexity is often delicate to maintain.

Motivated by this trade-off, we introduce a level-dependent intergrid strategy designed to balance efficiency and operator complexity.
Denote the standard full-weighting restriction (corresponding to bilinear interpolation) by (see \cite{trottenberg2000multigrid})
\begin{equation}\label{eq:biliniearIntergrid}
R_\mathrm{bilin} = \frac{1}{16}
\begin{bmatrix}
1 & 2 & 1 \\
2 & 4 & 2 \\
1 & 2 & 1
\end{bmatrix}
\end{equation}
and the bicubic restriction by (see \cite{holtz2007b})
\begin{equation}\label{eq:bicubicIntergrid}
R_\mathrm{bicub} = \frac{1}{256}
\begin{bmatrix}
1 & 4 & 6 & 4 & 1 \\
4 & 16 & 24 & 16 & 4 \\
6 & 24 & 36 & 24 & 6 \\
4 & 16 & 24 & 16 & 4 \\
1 & 4 & 6 & 4 & 1
\end{bmatrix}.
\end{equation}
The corresponding prolongation operators are $P_\mathrm{bilin} = 4R_\mathrm{bilin}^T$ and $P_\mathrm{bicub} = 4R_\mathrm{bicub}^T$.
We suggest the following level-dependent intergrid scheme, in which the restriction between coarse levels is \emph{not} a scaled transpose of the prolongation: for the intergrid between the first and second level, 
\begin{equation}\label{eq:level1Intergrid}
R_{1\to2} = R_\mathrm{bicub}, \quad P_{2\to1} = P_\mathrm{bicub}
\end{equation}
and between the $k$-th and $(k+1)$-st levels (where $k>1$)
\begin{equation}
\label{eq:levelDepIntergrid}
R_{k\to k+1} = R_\mathrm{bilin}, \quad P_{k+1\to k} = P_\mathrm{bicub}.
\end{equation}

As observed in \cite{chen2024matrix}, when using the bilinear intergrid \eqref{eq:biliniearIntergrid}, the computational stencil of the Galerkin operator remains $3\times 3$ in each level, while when using the bicubic intergrid \eqref{eq:bicubicIntergrid}, the stencil grows to $5\times 5$ in the second level and $7\times 7$ in the rest of the levels.
Using a mixed scheme --- bilinear restriction and bicubic interpolation --- limits the growth of the computational stencil, and it does not exceed $5\times 5$ at any level.
Yet, our experience shows that the multigrid convergence using mixed intergrid operators deteriorates compared to bicubic operators, for high-frequency Helmholtz problems. 
Hence, the need for a scheme that limits the operator complexity without hampering convergence gives rise to the level dependent scheme described in \eqref{eq:level1Intergrid} and \eqref{eq:levelDepIntergrid}.
The method builds on the understanding that the first level has the most significant influence on the entire cycle, and as we show later in Section \ref{sec:results}, investing more effort in the first level pays off.

The generalization of the suggested intergrid scheme to 3D is straightforward: the trilinear and tricubic intergrid operators are defined, similarly to the bilinear and bicubic intergrid operators, by the repeating Kronecker product
\begin{equation} \label{eq:trilin}
R_\mathrm{trilin} = \frac{1}{4^3}\begin{bmatrix}
1 & 2 & 1
\end{bmatrix}
\otimes 
\begin{bmatrix}
1 & 2 & 1
\end{bmatrix}
\otimes 
\begin{bmatrix}
1 & 2 & 1
\end{bmatrix}
\end{equation}
and
\begin{equation} \label{eq:tricub}
R_\mathrm{tricub} = \frac{1}{16^3}\begin{bmatrix}
1 & 4 & 6 & 4 & 1
\end{bmatrix}
\otimes 
\begin{bmatrix}
1 & 4 & 6 & 4 & 1
\end{bmatrix}
\otimes 
\begin{bmatrix}
1 & 4 & 6 & 4 & 1
\end{bmatrix}.
\end{equation}
The corresponding level-dependent intergrid scheme is defined similarly to \eqref{eq:level1Intergrid} and \eqref{eq:levelDepIntergrid}.
The 3D framework clearly inherits the compactness properties of the corresponding Galerkin coarse grid operators, namely, the coarse grid operators resulting from $H_c = RHP$ with level dependent $R$ and $P$, have at most $5\times 5 \times 5$ computational stencils, regardless of the number of levels.
Compared to $7\times 7 \times 7$ computational stencils resulted by tricubic intergrid operators, the save of computations is even more significant in 3D.

To measure the computational cost of each of the above mentioned intergrid schemes, we calculate its operator complexity as a function of the number of levels. 
Below, we use the notation $H_c^1=H$ for the finest grid operator and $H_c^l$ for the resulting Galerkin operator of the $l$-th level (given a certain intergrid scheme).
The operator complexity is calculated by
\begin{equation}\label{eq:operatorComplexity}
\frac{\sum_{l=1}^{m}(\text{nnz}(H_c^l))}{\text{nnz}(H)}
\end{equation}
where $m$ is the number of levels and ``nnz'' is a function that counts the nonzeros.
Table \ref{tab:operator_complexity} demonstrates the limited growth of the Galerkin operators for different choices of intergrid operators in 3D.

\begin{table}
\centering
\begin{tabular}{c|ccccc}
\hline
  \toprule
  \mc{5}{c}{Operator complexity}\\
  \midrule
 & 2-level & 3-level & 4-level  & 5-level \\
  \midrule
tricubic & 1.789 & 2.030 & 2.056 & 2.058 \\
mixed / lev-dep  & 1.789 & 1.886 & 1.898 & 1.899 \\
trilinear & 1.179 & 1.202 & 1.202 & 1.206 \\
  \bottomrule
 \end{tabular}
\caption{The operator complexity, using \eqref{eq:operatorComplexity}, of 2, 3, 4, and 5-level methods using GCA and cubic intergrid, for a $64\times 64\times 64$ cells grid in 3D.}
\label{tab:operator_complexity}
\end{table}

\subsection{Additive Vanka smoothing}\label{subsec:smoother}

Vanka relaxation was originally developed for saddle-point systems that arise in incompressible fluid flow simulations discretized on a staggered grid \cite{vanka1986blockFlow}. In this context, all unknowns at different locations within the same cell are gathered, and the resulting small submatrices can be relaxed at once with a proper weighting in the overlapping nodes, resulting in an additive Vanka method, or the submatrices can be relaxed sequentially, resulting in a multiplicative Vanka method.
Vanka relaxation methods can also be viewed as a particular class of overlapping domain decomposition methods, in which local solves are performed on small subdomains associated with degrees of freedom or grid patches.

In the context of multigrid, Vanka relaxation was shown to be a very efficient smoother for a wide range of applications. 
Additive Vanka smoothers are applied in \cite{farrell2021local} to the Stokes equation, and in \cite{he2022parameter} for linear elasticity problems.
Multiplicative Vanka smoothers appear in various applications, such as heat transfer, electromagnetic fields, incompressible fluid dynamics, and more.
Recently, an adaptation of Vanka smoothers for unstructured grid was suggested in \cite{nytko2025unstructured}
, by mapping the computations onto a structured grid.
Particularly, LFA for such smoothers is given in
\cite{sivaloganathan1991use, rodrigo2016local, maclachlan2011local, de2019robust, de2021two}, and particularly for the elastic Helmholtz equation in
\cite{treister2024hybrid}.
Despite being developed for systems of equations discretized on a staggered grid, Vanka smoothers can be used also for scalar equations discretized on nodal grids:
in \cite{greif2023closed}, different choices of patches for an additive Vanka smoother are applied to the Poisson equation, including a closed form analysis.

In this work we focus on additive Vanka smoothers, which can be applied in parallel and can be analyzed more simply.
We follow the formalism of \cite{greif2023closed} regarding additive Vanka smoothers for scalar equations, and suggest new patches that fit the Helmholtz equation \eqref{eq:acousitcHelm} with the discretization \eqref{eq:disc4thCompact}.
We equip these patches with special treatment of  near-boundary patches, which are not necessary for Poisson-like equations, yet enhance the convergence for the Helmholtz equation.

Let $X$ be the set of degrees of freedom (DOF) in the given grid, and let $X_i$ for $i=1,...,m$ be the (possibly intersecting) sets of DOFs in each patch, 
such that $X=\bigcup_{i=1}^m X_i$.
Let $V_i$ be an injection operator that maps a grid vector to the corresponding vector on the nodes in $X_i$. 
Then the $i$-th patch matrix is defined as
\begin{equation}
H_i = V_i H V_i^T
\end{equation}
with $H$ representing the discretized matrix of \eqref{eq:acousitcHelm}.
We update the current iterate $\bfu^{(k)}$ by restricting the residual, solving an error residual equation in each patch's space and updating the iterate using the corrected error. 
Namely,
\begin{equation}
H_i\bfe_i = V_i\left(\bfq - H\bfu^{(k)}\right)
\qquad \text{and} \qquad
\bfu^{(k+1)} = \bfu^{(k)} + W_i \bfe_i 
\end{equation}
where $W_i$ is a diagonal weighting matrix that spreads the weight of overlapping nodes between all the patches $X_i$ that include the same node. 
Typically, we take a scalar matrix $W = \frac{1}{N}I$ where $N$ is the size of the patch, and each node participates in $N$ patches.

Finally, an additive Vanka relaxation error can be interpreted as
\begin{equation}\label{eq:smoother}
S = I - w \left(\sum_{i=1}^m V_i^T W_i H_i^{-1} V_i \right) H
\end{equation}
where $w$ is a damping parameter.

\begin{figure}
\begin{center}
	\newcommand{\image}[1]{\includegraphics[width=0.185\linewidth]{#1}}
    \subfigure[\footnotesize Element patch]{\image{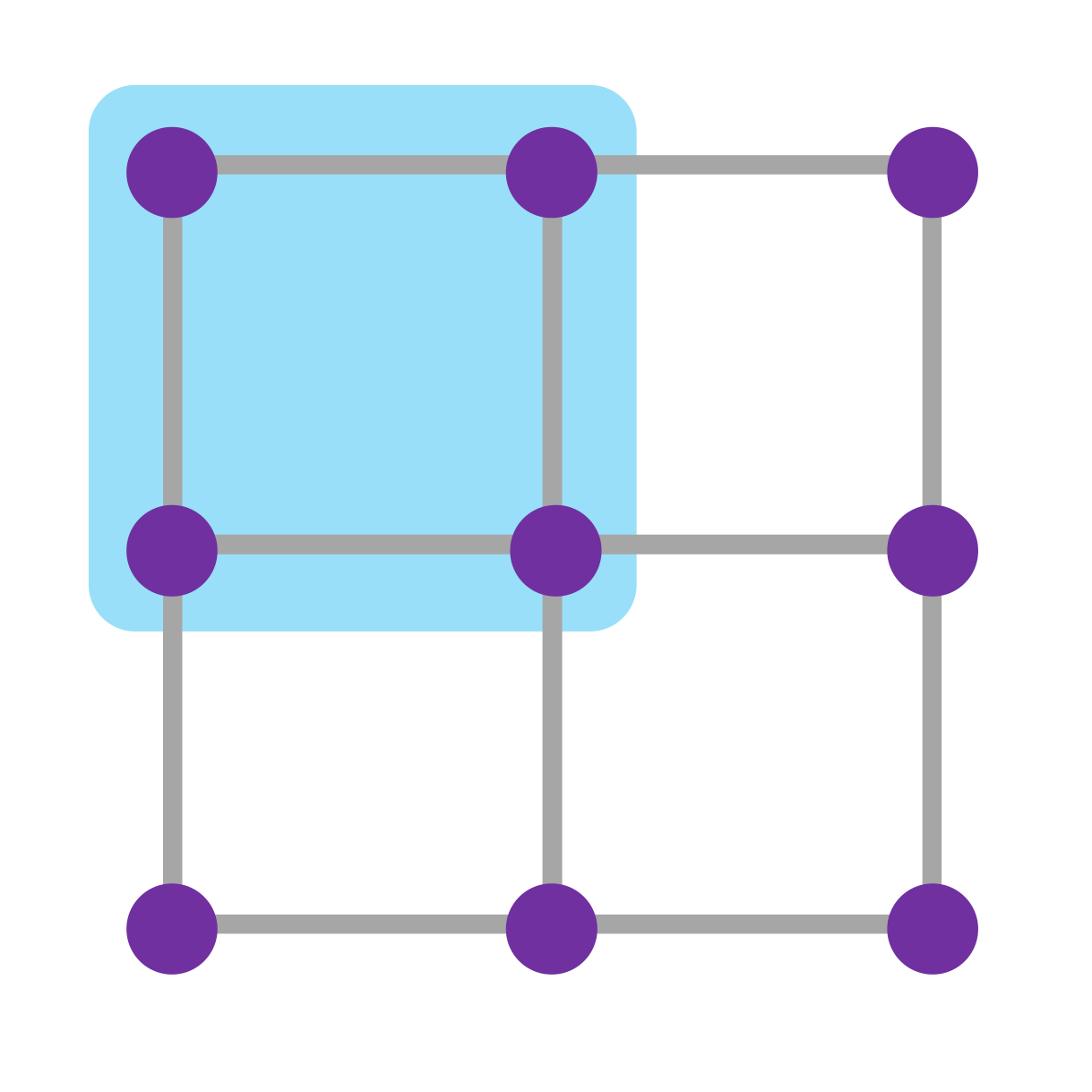}\label{fig:4patch}} 
    \hspace{20pt}
    \subfigure[\footnotesize Plus patch]{\image{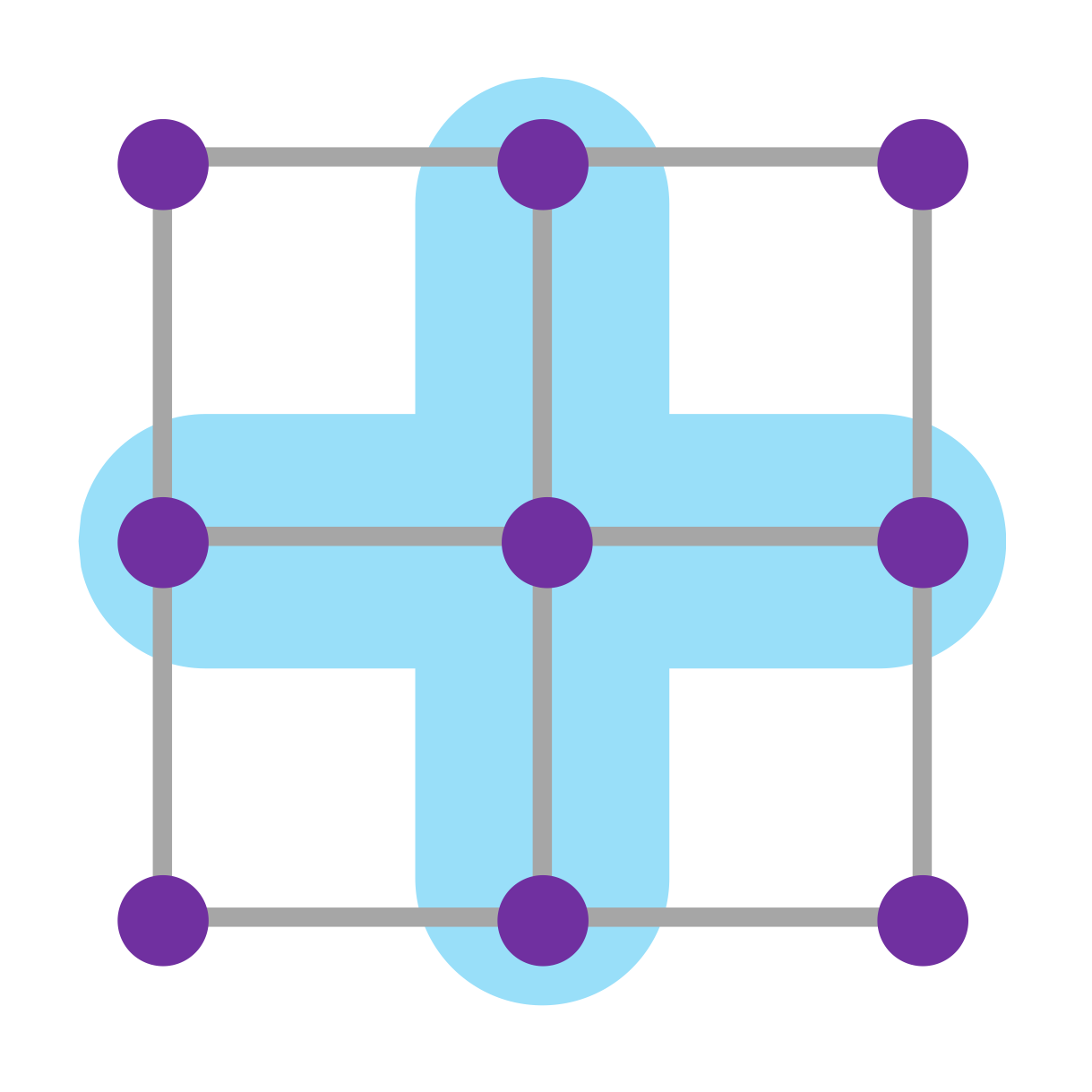}\label{fig:5patch}}
    \hspace{20pt}
    \subfigure[\footnotesize RB patch]{\image{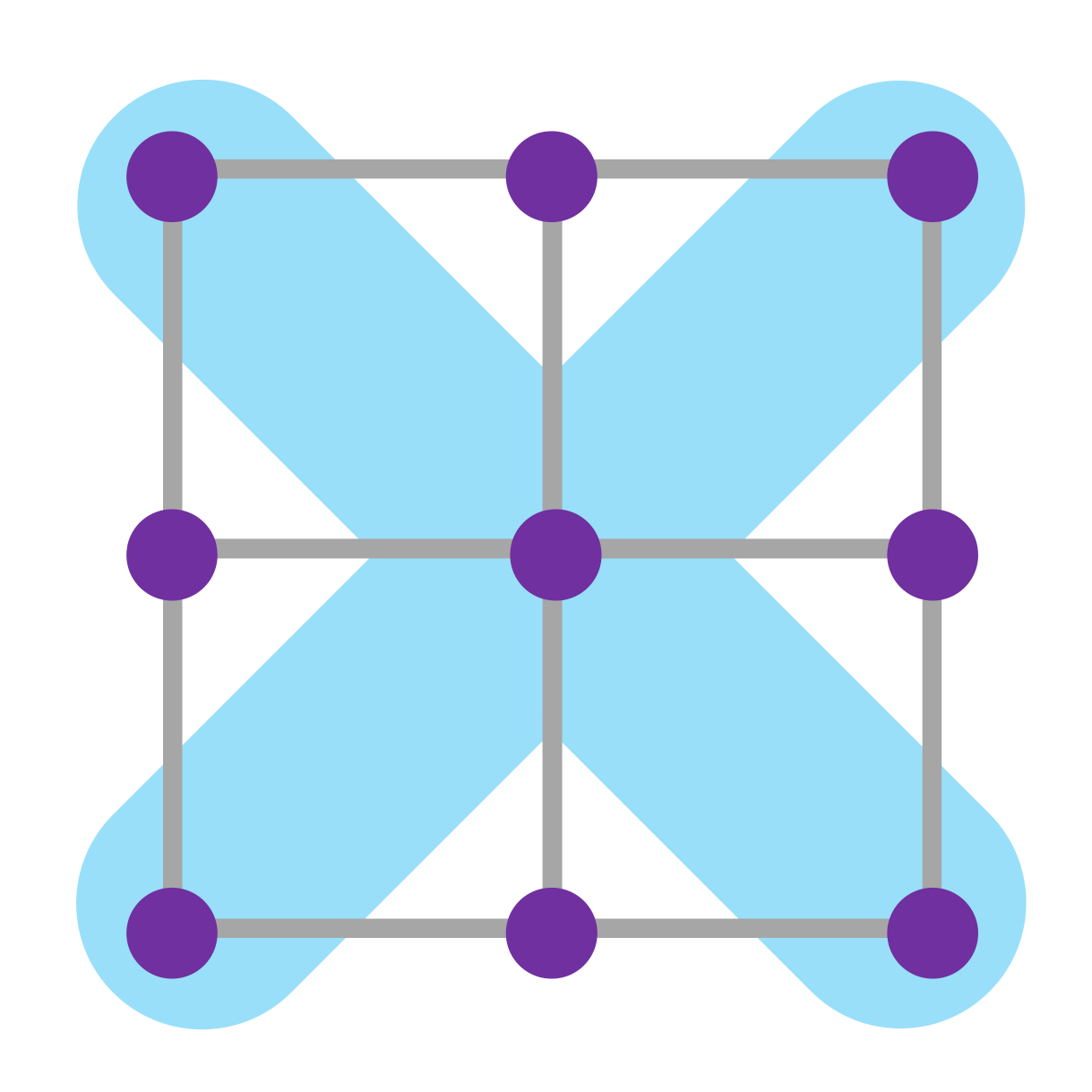}\label{fig:skew5patch}}
    \hspace{20pt}
    \subfigure[\footnotesize Full patch]{\image{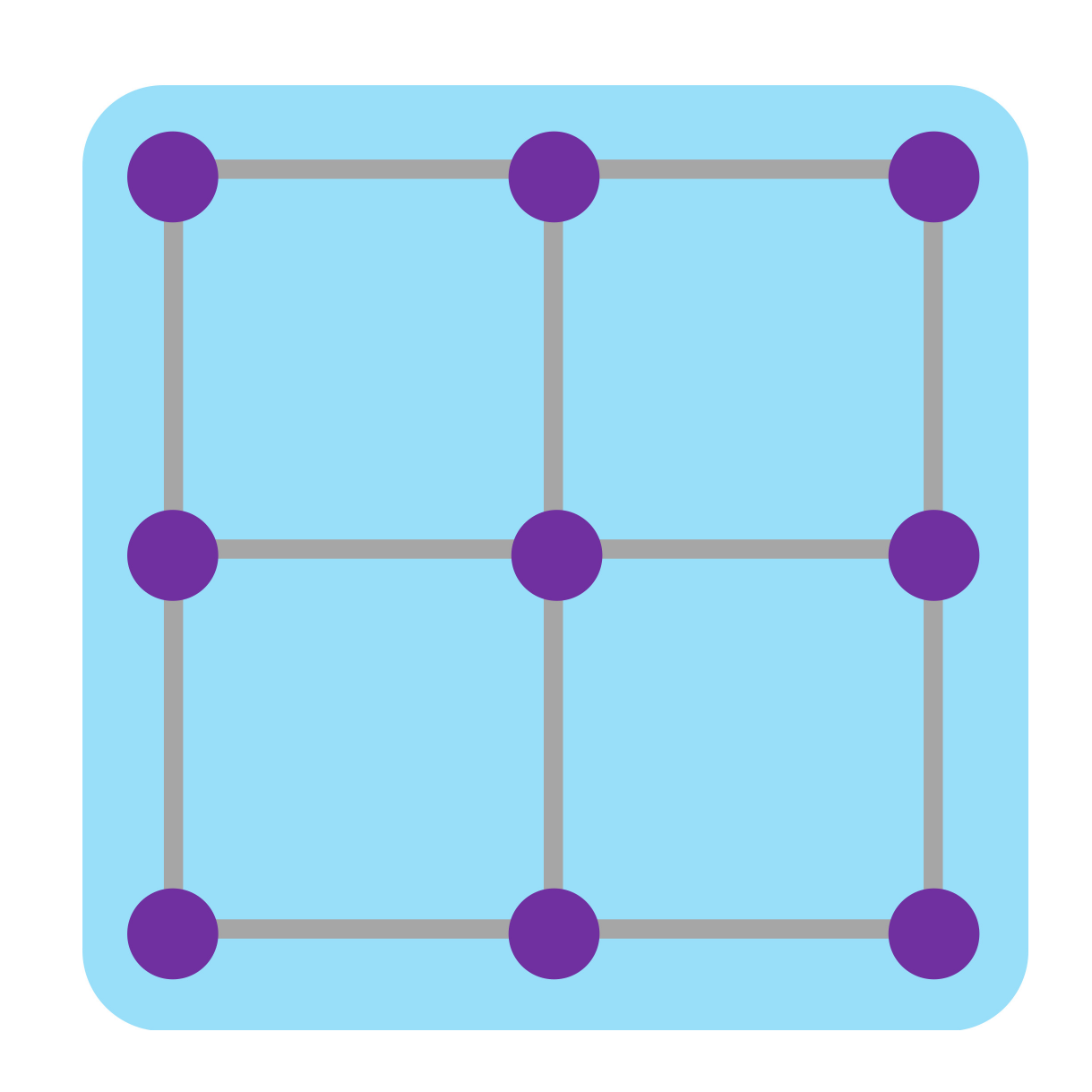}\label{fig:9patch}}\\
\end{center}
\caption{Patches for additive Vanka smoother for a 2D nodal discretization.
}\label{fig:patches}
\end{figure}

\begin{figure}
\begin{center}	\newcommand{\image}[1]{\includegraphics[width=0.23\linewidth]{#1}}
    \subfigure[\footnotesize Element patch]{\image{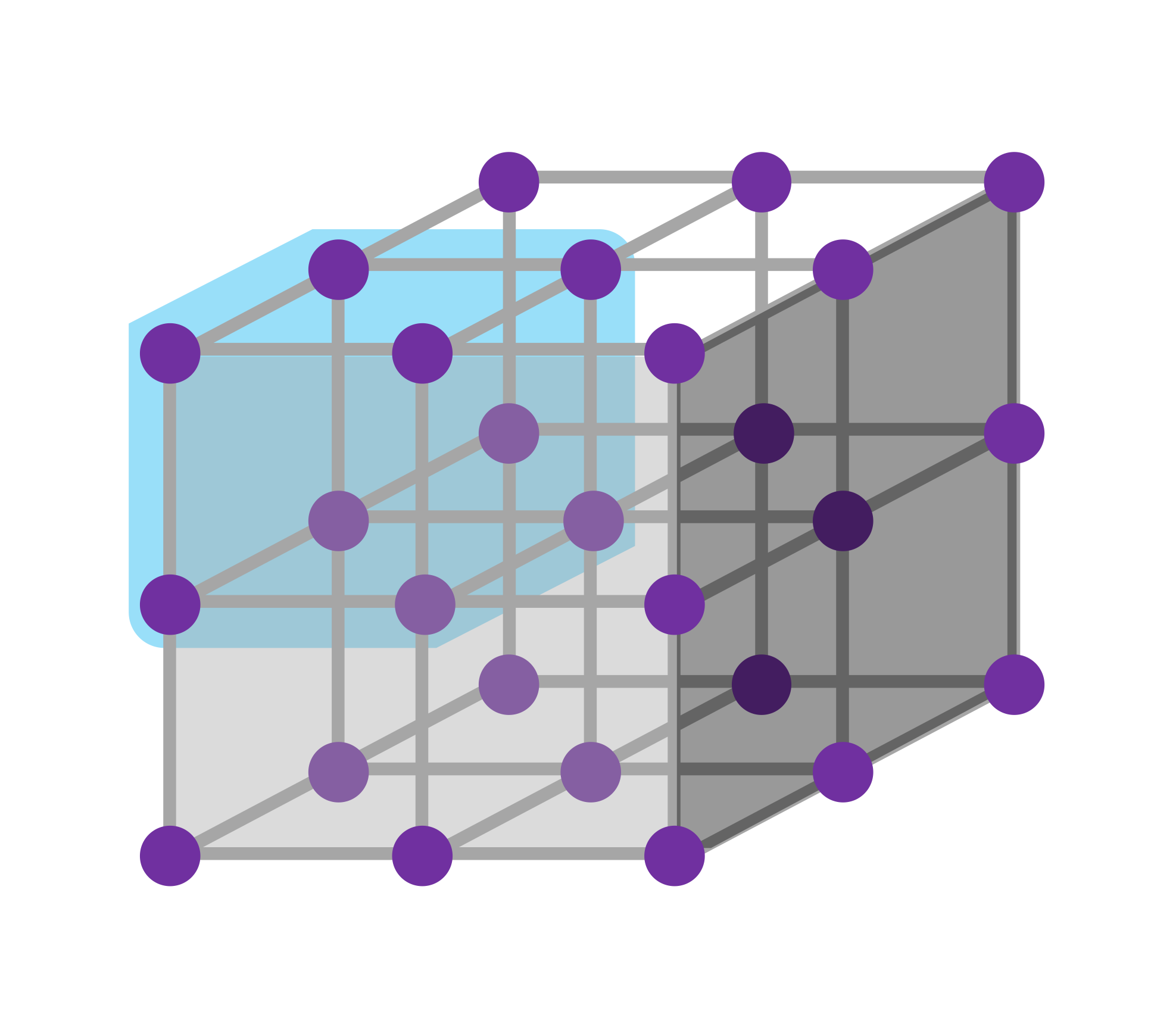}\label{fig:8patch}} 
    \subfigure[\footnotesize Plus patch]{\image{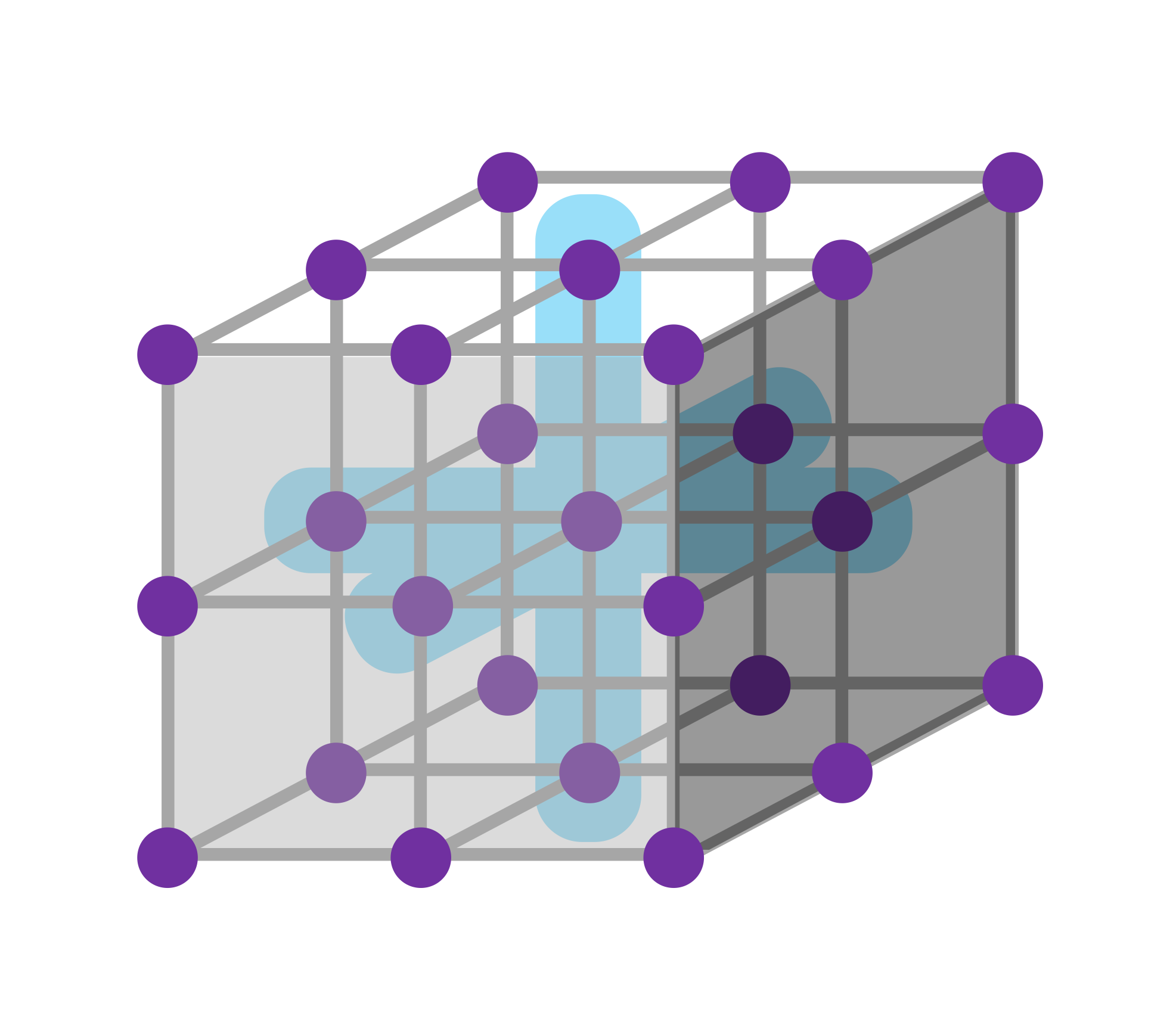}\label{fig:7patch}}
    \subfigure[\footnotesize RB patch]{\image{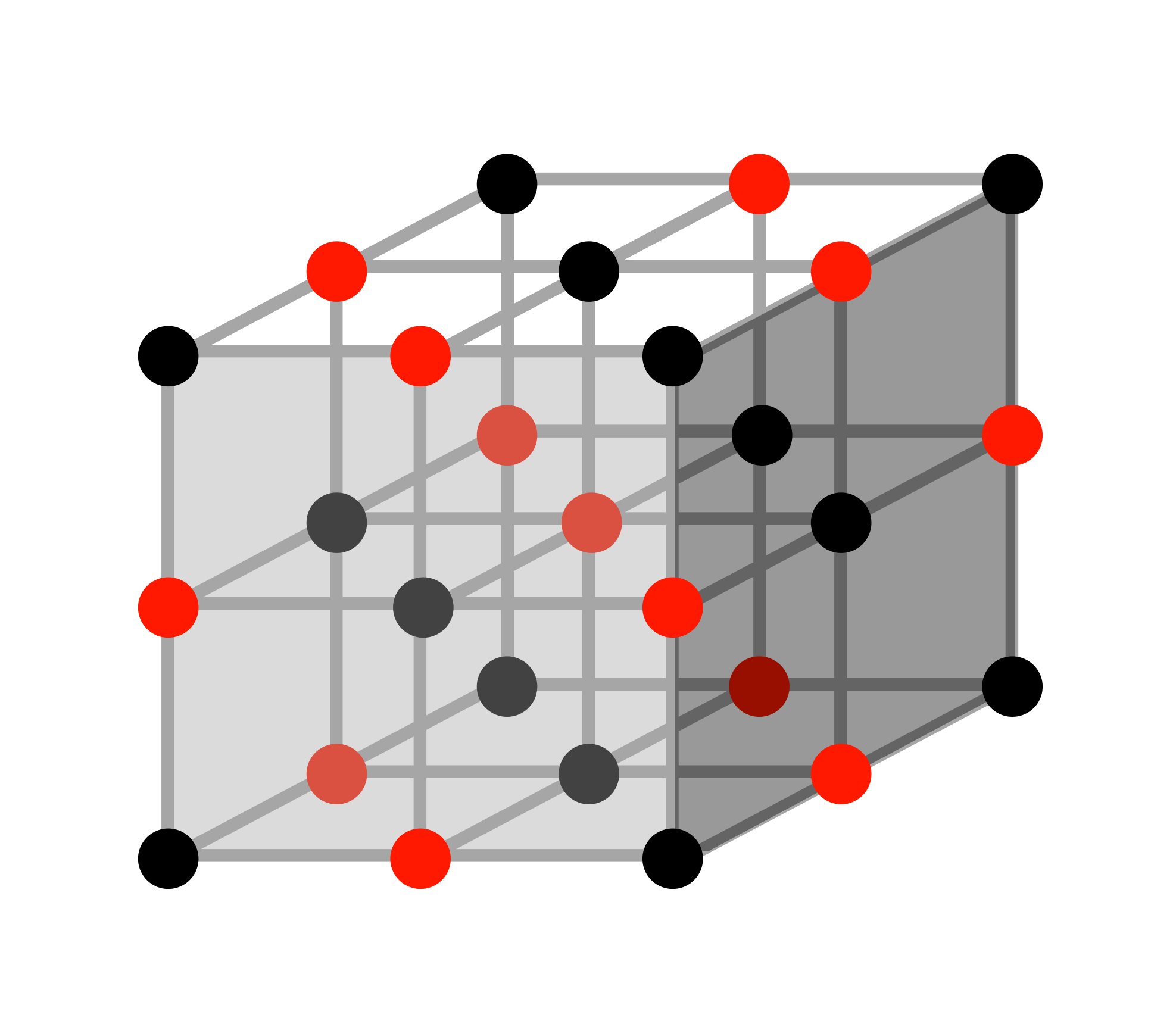}\label{fig:13patch}}
    \subfigure[\footnotesize Full patch]{\image{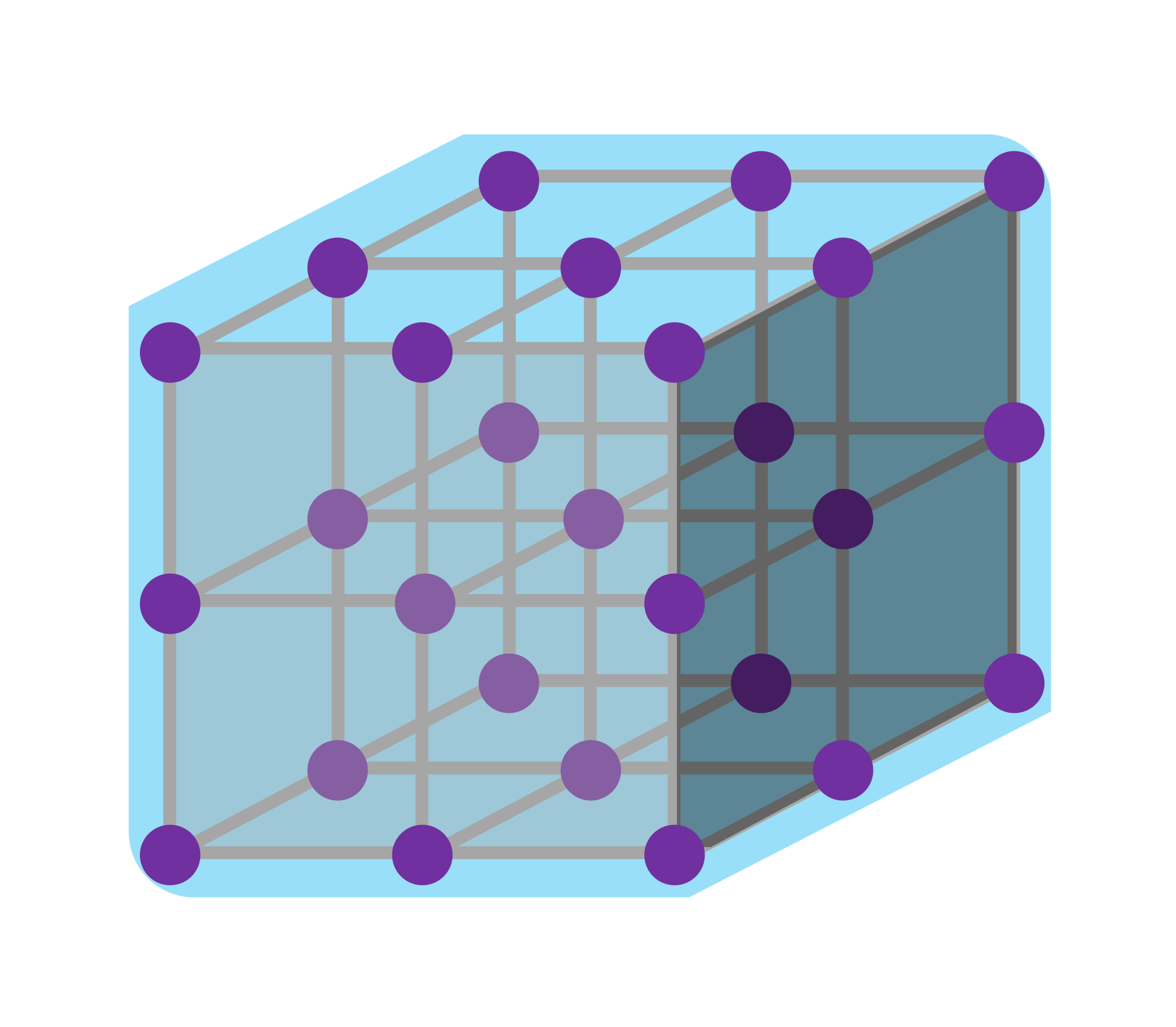}\label{fig:27patch}}
\\
\end{center}
\caption{Patches for additive Vanka smoother for a 3D nodal discretization.
}\label{fig:patches3D}
\end{figure}

Figs. \ref{fig:patches} and \ref{fig:patches3D} 
describe the following Vanka patches in 2D and 3D, respectively.
\begin{itemize}
\item \textbf{Element patch:} includes all the vertices in one discretization cell, which are 4 vertices in 2D, see Fig. \ref{fig:4patch}, and 8 vertices in 3D, see Fig. \ref{fig:8patch}.
\item \textbf{Plus patch:} includes all the vertices that communicate with the central vertex via a 2nd order discretization. That is, plus shaped 5 vertices in 2D, see Fig. \ref{fig:5patch}, and 7 vertices in 3D, see Fig. \ref{fig:7patch}.
\item \textbf{RB patch:} includes all the vertices of the same color as the central vertex in a red-black coloring.
In 2D, it comprises an X-shaped patch with 5 vertices, see Fig. \ref{fig:skew5patch}, and in 3D a patch of 13 vertices, see Fig. \ref{fig:13patch}.
\item \textbf{Full patch:} includes all the vertices that communicate with the central vertex via a compact 4-th order stencil, such as the stencils that we use, \eqref{eq:disc4thCompact} in 2D and \eqref{eq:disc4thCompact3D} in 3D.
The 2D full patch, comprised of 9 vertices, is depicted in Fig. \ref{fig:9patch}, and its 3D analogue, comprised of 27 vertices, is in Fig. \ref{fig:27patch}.
\end{itemize} 
The element patch and the plus patch were investigated in \cite{greif2023closed} in the context of a second-order discretization, and the plus patch is named there ``vertex-wise patch''. 
The wider discretization we use here raises the need for wider patches, such as the Full patch, and its sub-patch, the RB patch.
In Fig. \ref{fig:stencils} we depict the pattern of the Galerkin stencils of the Laplacian in 2D, which inspires the definition of the RB patch.

\begin{figure}
\begin{center}
	\newcommand{\image}[1]{\includegraphics[width=0.22\linewidth]{#1}}
    \subfigure[\footnotesize Fine level]{\image{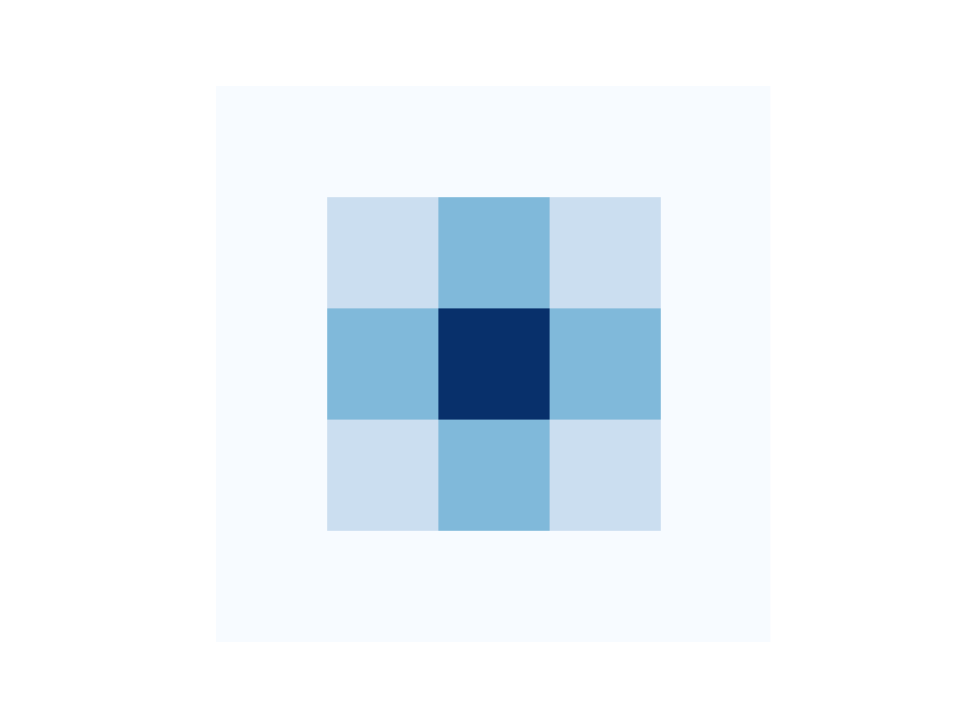}\label{fig:stencil_level_1}} 
    \subfigure[\footnotesize Second level]{\image{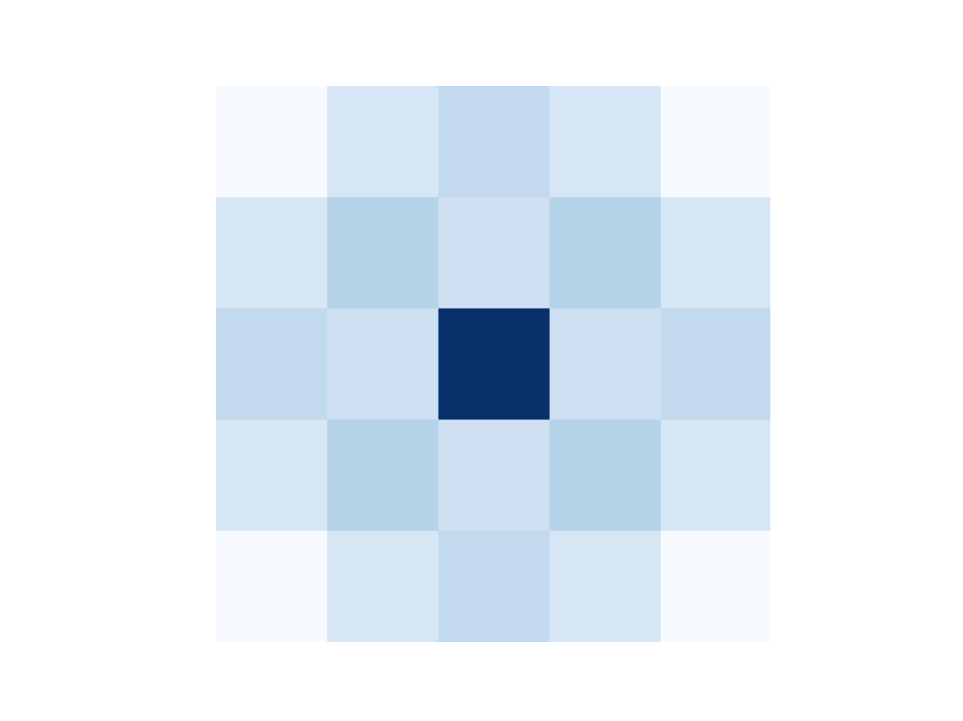}\label{fig:stencil_level_2}}
    \subfigure[\footnotesize Third level]{\image{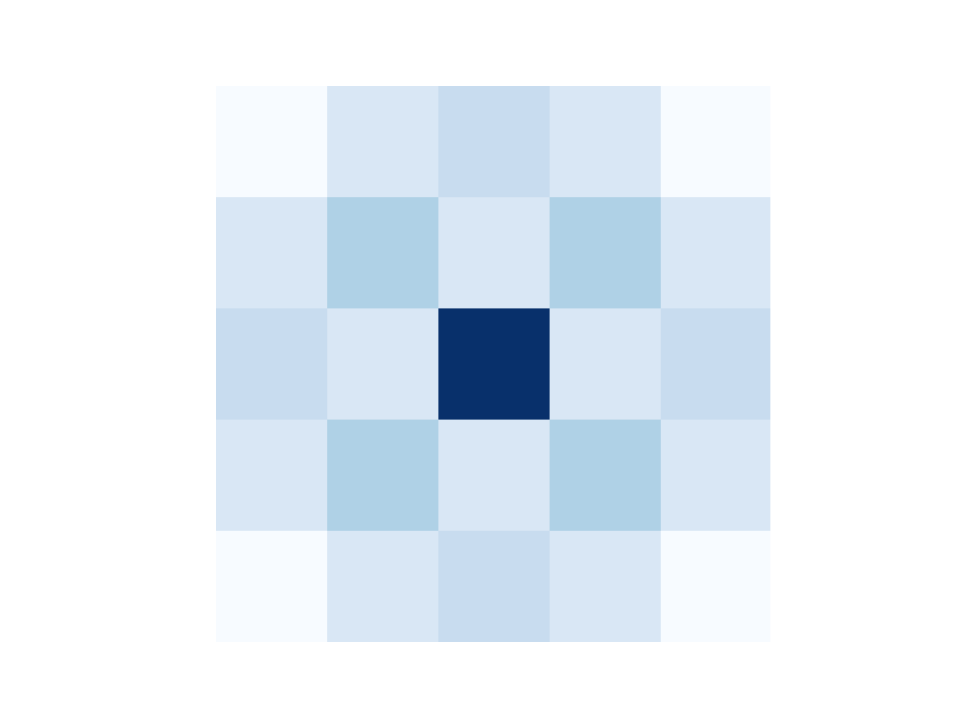}\label{fig:stencil_level_3}}
    \subfigure[\footnotesize Fourth level]{\image{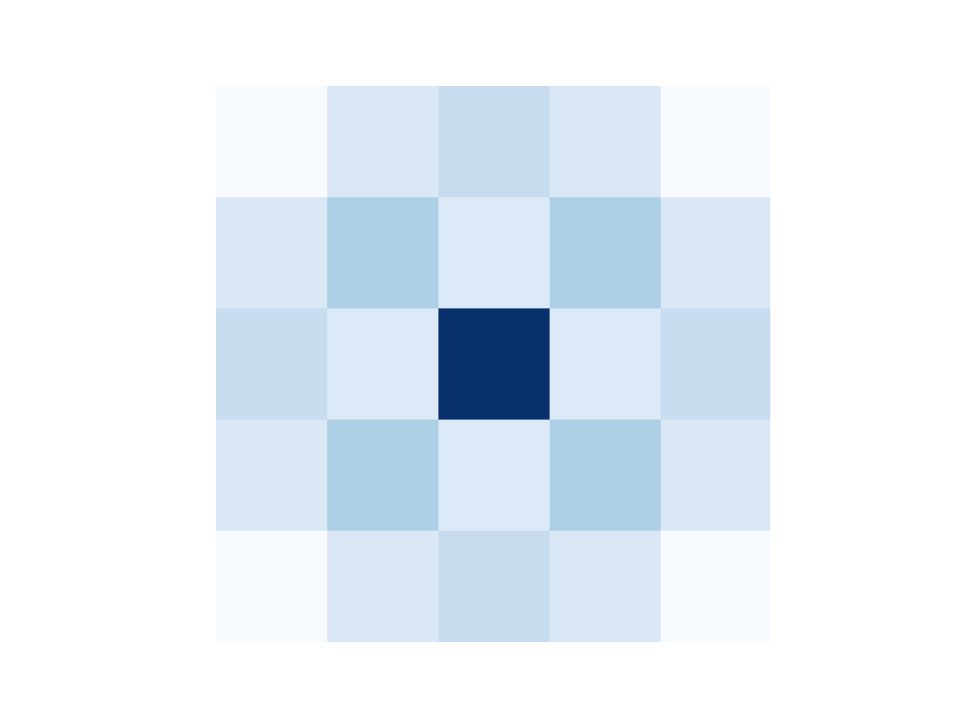}\label{fig:stencil_level_4}}\\
\end{center}
\caption{Magnitudes of the 2D Laplacian's stencil coefficients for different levels. 
The fine stencil is the compact 4-th order discretization described in \eqref{eq:disc4thCompact}, and the other stencils are accepted by Galerkin coarse approximation with the level-dependent intergrid scheme described in \eqref{eq:levelDepIntergrid}.
}\label{fig:stencils}
\end{figure}

\subsection{Treatment of boundaries}

Efficient implementation of additive Vanka smoothers requires special treatment of boundary, especially for indefinite problems such as the Helmholtz equation.
The special treatment we suggest includes two stages:
adding smaller patches on the boundary, and recalculating the weighting matrices $W_i$ from \eqref{eq:smoother} (rather than using scalar matrices).
Without these adaptations, the method resembles point-wise Jacobi with an extremely small damping parameter near the boundaries.
For (weakly) diagonally dominant problems, such as the Poisson equation, where Jacobi is a stable smoother (and might even converge as a solver), no further treatment is needed.
However, we observed in our experiments that for high-frequency Helmholtz problems, these two stages of treatment are necessary to meet the theoretical predictions of convergence rate.

In the first stage, we add smaller boundary patches that are the intersection of an entire patch with the domain.
We only do it for the wide patches, that is, all the patches except for the element patch.
For the case of the 2D plus patch, e.g., it results in 3-point patches for each corner and triangle-shaped 4-point patches along the edges of the domain, see Fig. \ref{fig:boundary_5patch}. 
We note that similar patches were constructed  in \cite{claus2021nonoverlapping} in the context of non-overlapping smoothers for the Stokes equation.
For the 2D RB patch it results in 2-point patches for each corner and wedge-shaped 3-point patches along the edges of the domain, see Fig. \ref{fig:boundary_skew5patch}, and for the 2D Full patch, a similar reasoning leads to 4-point patches for each corner and rectangular 6-point patches along the edges, see Fig. \ref{fig:boundary_9patch}.

In the second stage, we recalculate the weights for the weighting matrices $W_i$.
First, for each node we count the number of patches that include it. 
This number is written inside each node in Fig. \ref{fig:boundary_patches}. 
Second, we define the weight of this node as the inverse of this number.
Finally, we take the $(j,j)$-th element in $W_i$ to be the weight of the $j$-th node in the $i$-th patch.
Clearly, for patches that do not intersect with the boundaries, this process leads to the equal weighting $W=\frac{1}{N}I$, where $N$ is the size of the patch.
In 3D, we treat the boundary in an analogous manner.

\begin{figure}
\begin{center}
	\newcommand{\image}[1]{\includegraphics[width=0.186\linewidth]{#1}}
    \subfigure[\footnotesize Element patch]{\image{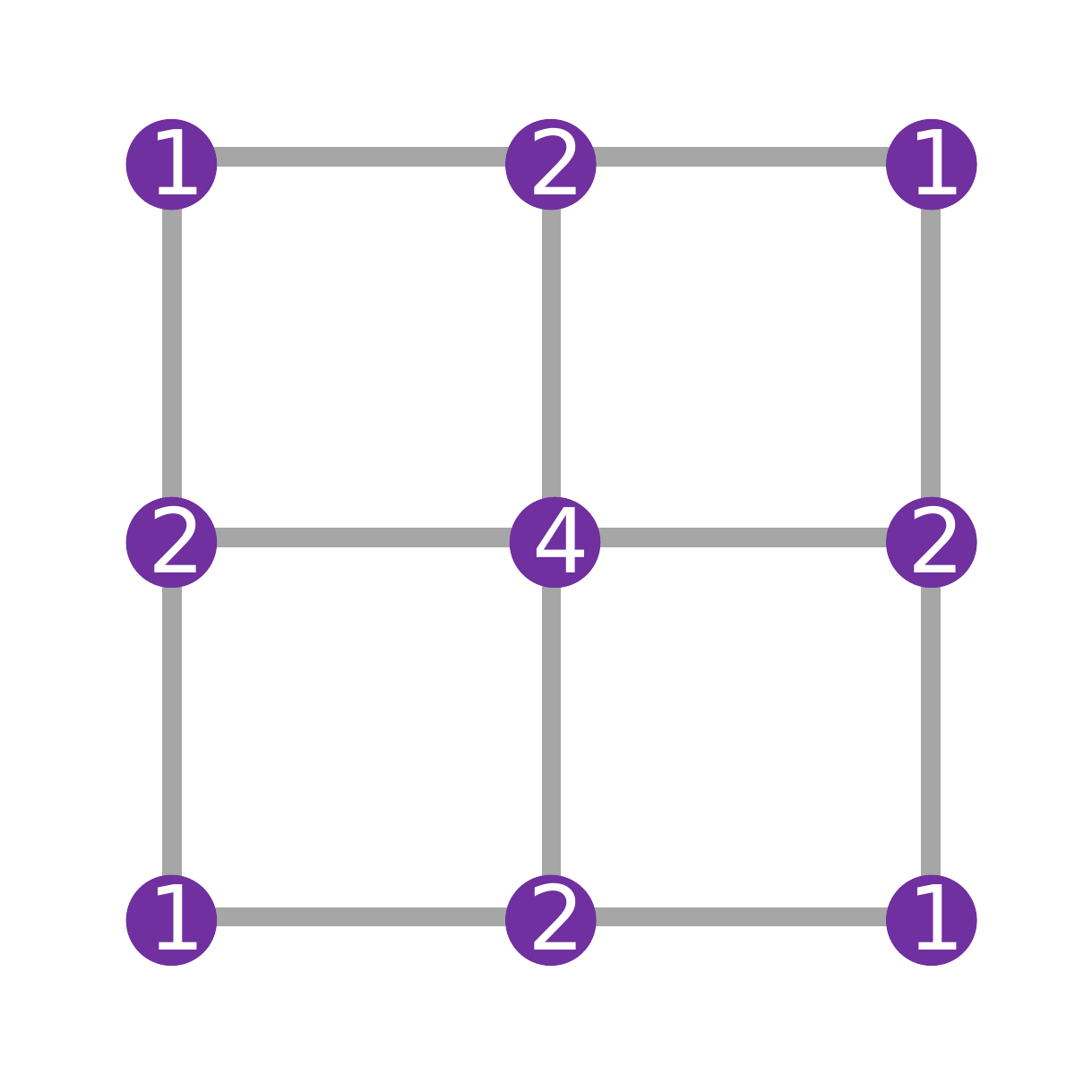}\label{fig:boundary_4patch}} 
    \hspace{20pt}
    \subfigure[\footnotesize Plus patch]{\image{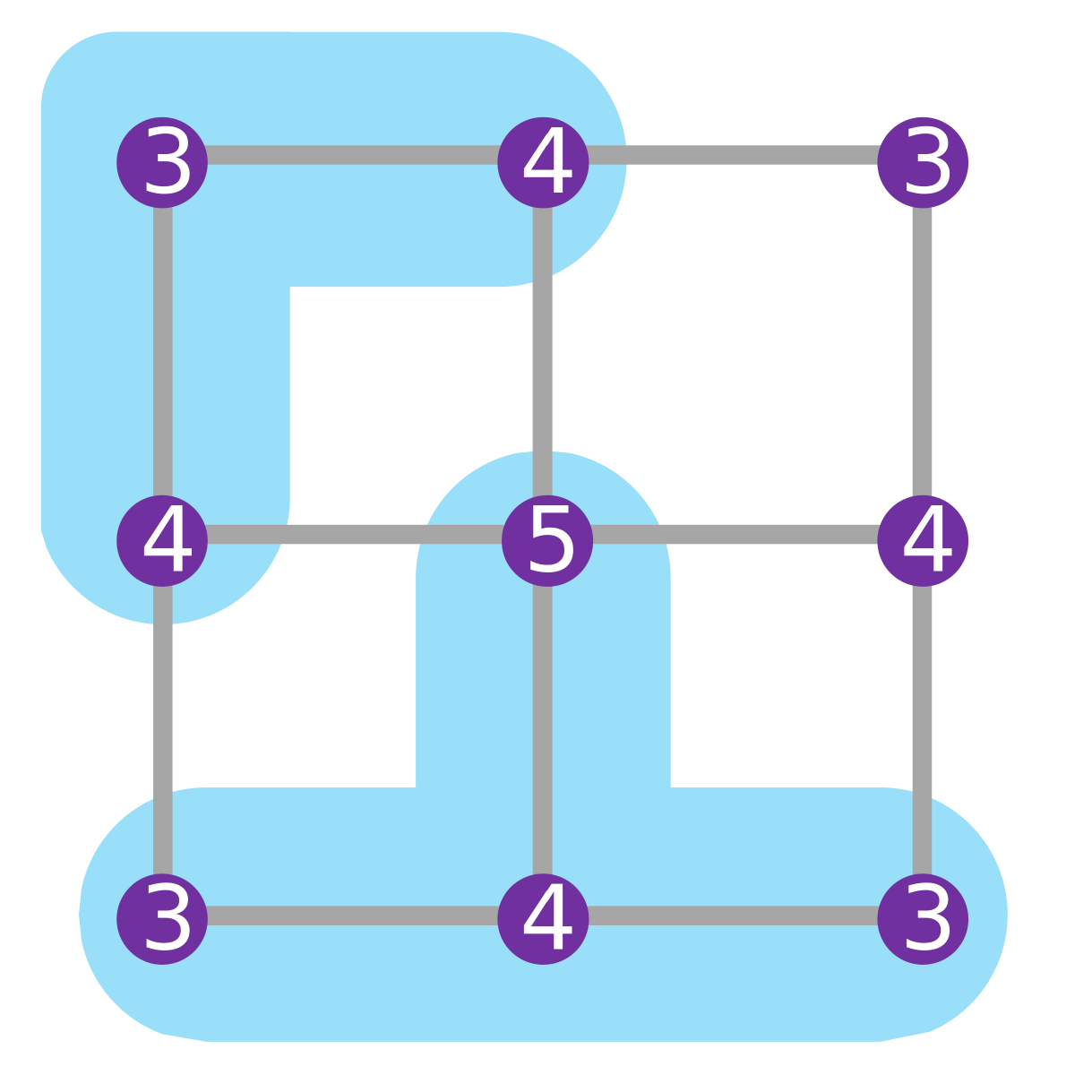}\label{fig:boundary_5patch}}
    \hspace{20pt}
    \subfigure[\footnotesize RB patch]{\image{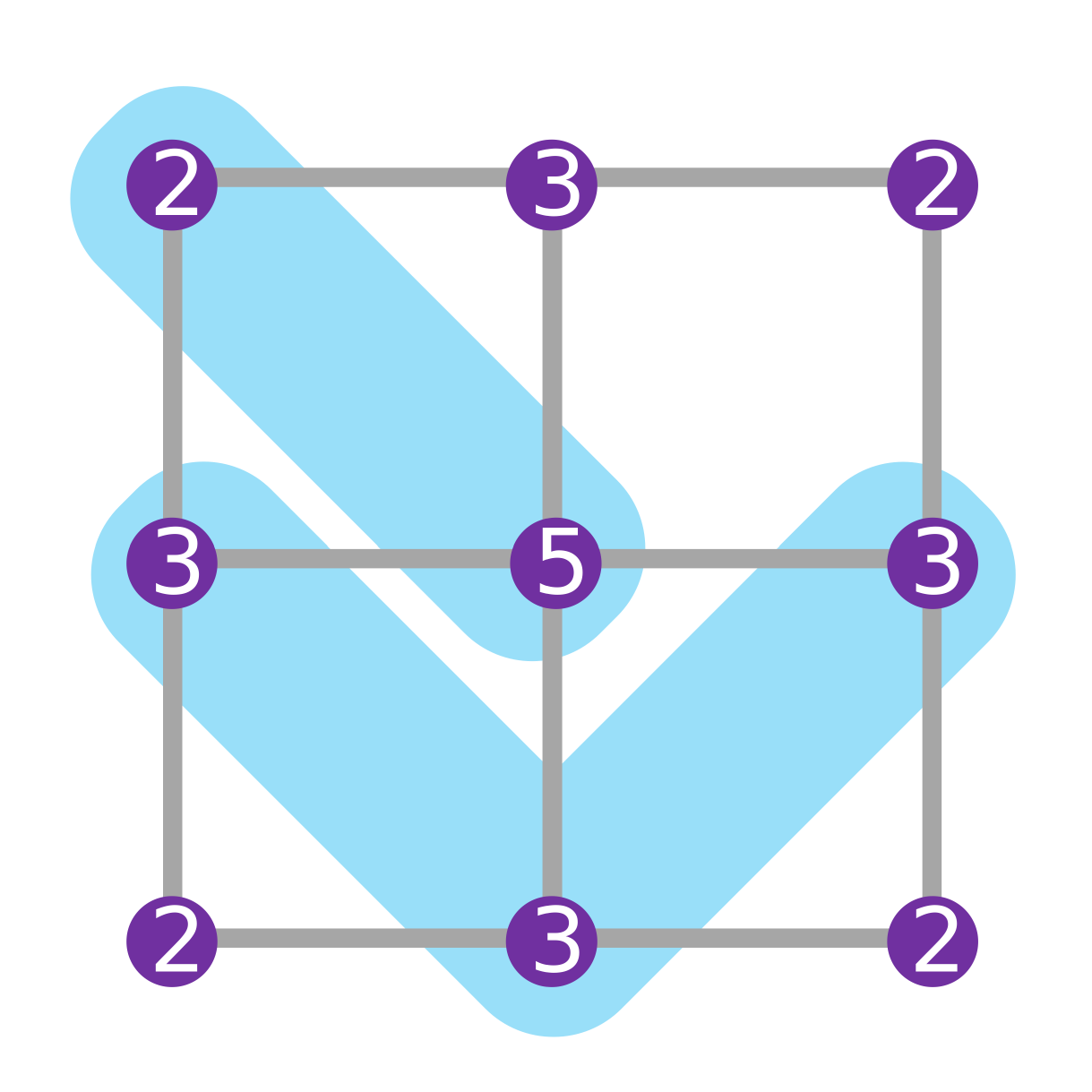}\label{fig:boundary_skew5patch}}
    \hspace{20pt}
    \subfigure[\footnotesize Full patch]{\image{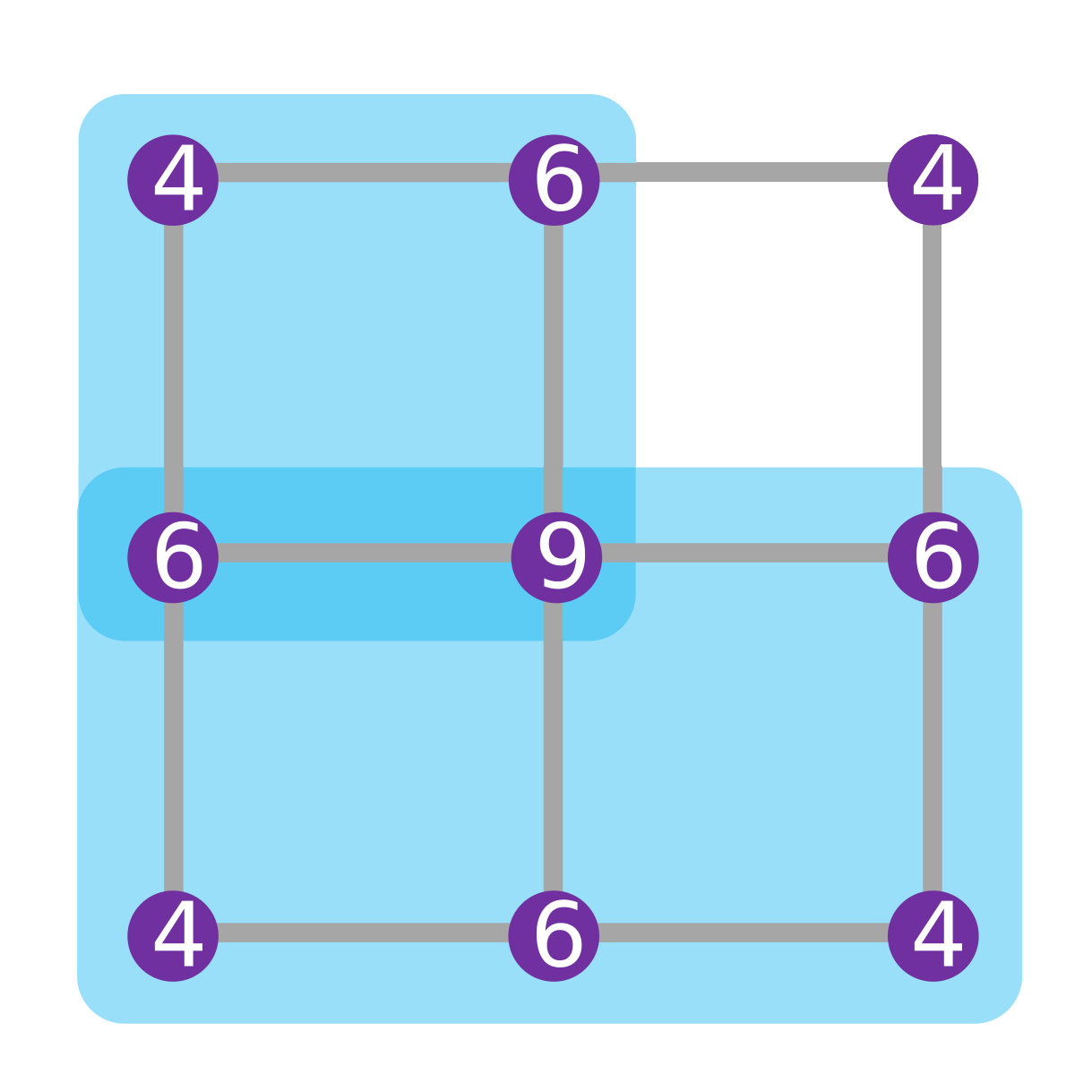}\label{fig:boundary_9patch}}\\
\end{center}
\caption{Additional boundary patches and boundary weights for additive Vanka smoother for a 2D nodal discretization.
}\label{fig:boundary_patches}
\end{figure}

\section{Local Fourier analysis}\label{sec:lfa}

In this section we describe the LFA calculation needed to analyze the described method. 
We first hold smoothing analysis for the additive Vanka patches described in Fig. \ref{fig:patches}, and then briefly describe the resulting two-grid analysis given the level-dependent intergrid operators.

\subsection{Smoothing LFA for additive Vanka}

Following the derivation in \cite{greif2023closed}, we calculate the symbol of the smoother as
\begin{equation}\label{eq:smoother_symbol}
\widetilde{S}(\theta) = 1 - w \left( V^T W \Phi(\theta)^H H_i^{-1} \Phi(\theta) V \right) \widetilde{H}(\theta)
\end{equation}
where $\widetilde{H}$ is the (scalar) symbol of $H$ from \eqref{eq:disc4thCompact}, $w$ is the damping, $V$ is an $N\times 1$ ones vector where $N$ is the number of points in the patch, $W$ is a diagonal weighting matrix that determines the relative weight of each overlapping node and $\Phi$ is a diagonal matrix with the Fourier modes corresponding to the patch locations on its diagonal.
Note that the resulting $\tilde{S}$ is a scalar symbol.

The symbol $\widetilde{H}$ depends only on the discretization.
Denote by $a=\frac{10}{3h^2} - \frac{2\sigma}{3}$, $b=-\frac{2}{3h^2} - \frac{\sigma}{12}$ and $c=-\frac{1}{6h^2}$ the discretization constants from the stencil \eqref{eq:disc4thCompact}, where $\sigma = \kappa^2\omega^2(1-\gamma\im)$. 
Then regardless of the patch, 
\begin{equation}\label{eq:Hsymbol}
\widetilde{H} = a + 2b \cos(\theta_1) + 2b \cos(\theta_2) + 4c \cos(\theta_1)\cos(\theta_2).
\end{equation}
For the element patch from Fig. \ref{fig:4patch}, we get 
\begin{equation}\label{eq:4lfa}
\Phi = diag\left[1,e^{\im \theta_1}, e^{\im\theta_2}, e^{\im (\theta_1+\theta_2)} \right]
\quad \text{and} \quad
H_i =
\begin{bmatrix}
a & b & b & c \\
b & a & c & b \\
b & c & a & b \\
c & b & b & a
\end{bmatrix}.
\end{equation}
for the plus patch from Fig. \ref{fig:5patch} we get
\begin{equation}\label{eq:5lfa}
\Phi = diag\left[e^{-\im \theta_2}, e^{-\im\theta_1},1, e^{\im \theta_1}, e^{\im\theta_2} \right]
\quad \text{and} \quad
H_i =
\begin{bmatrix}
a & c & b & c &   \\
c & a & b &   & c \\
b & b & a & b & b \\
c &   & b & a & c \\
  & c & b & c & a
\end{bmatrix},
\end{equation}
and for the RB patch from Fig. \ref{fig:skew5patch} we get
\begin{equation}\label{eq:skew5lfa}
\Phi = diag\left[e^{-\im(\theta_1 + \theta_2)}, e^{-\im(\theta_2-\theta_1)},1, e^{\im(\theta_2 - \theta_1)}, e^{\im(\theta_1 + \theta_2)} \right]
\end{equation}
and
\begin{equation}
H_i =
\begin{bmatrix}
a &  & c &  &   \\
 & a & c &   &  \\
c & c & a & c & c \\
 &   & c & a &  \\
  &  & c &  & a
\end{bmatrix}.
\end{equation}

We estimate the smoothing factor $\mu_\mathrm{loc}$ from Definition \ref{def_mu_loc} by sampling the symbol $\widetilde{S}(\theta)$ in $256\times256$ values of $\theta\in\left[-\frac{\pi}{2},\frac{3\pi}{2}\right]^2$, and maximizing numerically the resulting amplification factor over $T^\mathrm{high}$.
For the sake of comparison, the symbol for the Jacobi smoother is calculated as a $1\times 1$ Vanka patch, which is mathematically equivalent to the known symbol in the literature, see \cite{trottenberg2000multigrid}.
We note that the symbol of each smoother depends not only on the patch, but also on the underlying discretization stencil.

\subsection{Two-grid LFA}

To perform two-grid LFA, we calculate the symbol matrix $\widetilde{TG}$ from \eqref{eq:symbolTG}.
To calculate the symbol $\widetilde{S}(\theta)$, we evaluate \eqref{eq:smoother_symbol} for each of the four harmonics.
To calculate the operator's symbol $\widetilde{H}$ we use \eqref{eq:Hsymbol}.
To calculate the intergrid's vector of symbols, we first calculate the scalar symbol that corresponds to each intergrid stencil, and get
\begin{equation}
\widetilde{R}(\theta) = \frac{1}{4}\left(1+\cos(\theta_1)\right)\left(1+\cos(\theta_2)\right)
\end{equation}
for bilinear intergrid \eqref{eq:biliniearIntergrid} (see \cite{trottenberg2000multigrid}), and 
\begin{equation}
\widetilde{R}(\theta) =  \frac{1}{64}\left(3+4\cos(\theta_1)+\cos(2\theta_1)\right)\left(3+4\cos(\theta_2)+\cos(2\theta_2)\right).
\end{equation}
for bicubic intergrid \eqref{eq:bicubicIntergrid}.
Finally, to calculate the (scalar) symbol $\widetilde{H}_c$, given that $H_c=R H P$, we multiply the corresponding symbol matrices.
After calculating $\widetilde{TG}$, we sample it and maximize its value numerically in $T^\mathrm{low}$ to estimate the two-grid factor $\rho_\mathrm{loc}$ for a given  damping parameter $w$.
Using high-order intergrid enables a convergent two-grid cycle with no added shift, and hence we analyze $\rho_\mathrm{loc}$ for a non-attenuated system.

\section{Numerical results}\label{sec:results}

In this section we give numerical results.
Throughout the experiments, we solve the acoustic Helmholtz equation \eqref{eq:acousitcHelm} discretized on a rectangular domain.
We locate a point source in the center of the domain, and solve the problem up to a tolerance of relative residual $<10^{-6}$, with a zero initial guess.
We impose ABC as a layer of 20 grid cells of increasing attenuation towards the boundary.
We assume no physical attenuation, $\gamma = 0,$ in the rest of the domain.
In all of the experiments, we take a frequency that corresponds to $10$ grid points per wavelength. 

Our code is written in the {\tt Julia} language \cite{Julia}, and we include it as a part of the {\tt jInv.jl} package \cite{ruthotto2017jinv}.Our tests were computed on a dual-core laptop with 32 GB RAM, running Windows 11.
Our code is available online at 
\begin{center}
\tt https://github.com/RacheliYovel/VankaHelmholtz.git.
\end{center}

\subsection{LFA predictions vs. actual performance}

First, we compare the expected convergence by two-grid LFA to the convergence in practice.
We calculate $\rho_\mathrm{loc}$ as explained in Section \ref{sec:lfa}.
We compare it to the convergence factor calculated after $k$ iterations by
\begin{equation}\label{eq:cf}
c_f = \left(
\frac{\|\bfr^{(k)}\|}{\|\bfr^{(0)}\|} \right)^{1/k}
\end{equation}
where the calculation begins after a warm-up of 5 iterations.
We solve the 2D acoustic Helmholtz equation \eqref{eq:acousitcHelm} in a square dimensionless domain $\Omega=[0,1]^2,$ using a 2-level multigrid cycle with bicubic intergrid from equation \eqref{eq:bicubicIntergrid} and no added shift.
Note that the lack of shift enables the use of multigrid as a solver, and not only as a preconditioner, though this property holds only for 2-level cycles.
In order to demonstrate asymptotic behavior, we take the number of iterations $k$ that solve the problem up to a relative residual of $10^{-9}$.

\begin{figure}
\begin{center}
	\newcommand{\image}[1]{\includegraphics[width=0.3\linewidth]{#1}}
    \subfigure[\footnotesize Element patch]{\image{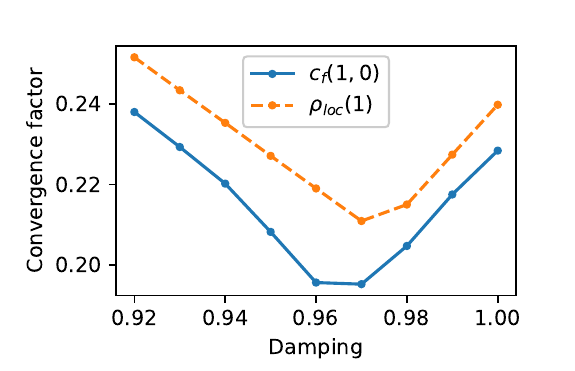}\label{fig:damping_4patch}} 
        \hspace{10pt}
    \subfigure[\footnotesize Plus patch]{\image{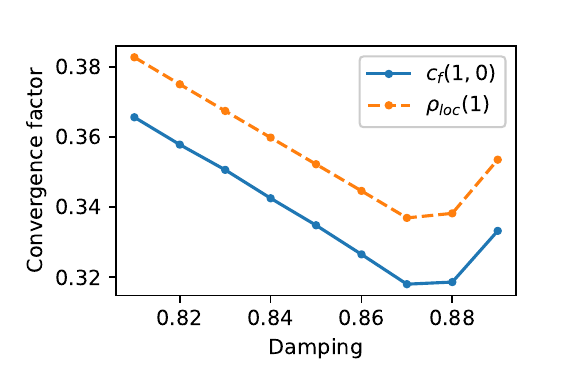}\label{fig:damping_5patch}}
    \hspace{10pt}
    \subfigure[\footnotesize RB patch]{\image{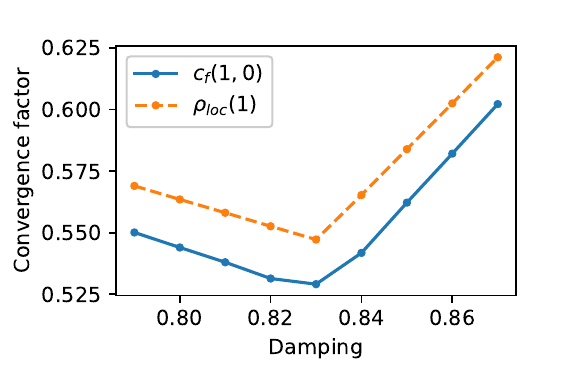}\label{fig:damping_skew5patch}}
\\
\end{center}
\caption{LFA two-grid factor and convergence factor in practice for different damping parameters $w$ for each patch in 2D.
}\label{fig:damping}
\end{figure}

In Fig. \ref{fig:damping} we compare $\rho_\mathrm{loc}$ and $c_f$ for different damping values, for homogeneous media with grid size $n_\mathrm{cells}=256\times256$, in order to find the optimal damping parameters for each patch from Fig. \ref{fig:patches}. 
We do not analyze here the Full patch from Fig. \ref{fig:9patch}, as the computational effort it requires per iteration is significantly higher: it has 9 points per patch, compared to 4--5 in the rest of the patches, and based on our experience the extra work does not pay off.
The small differences between $\rho_\mathrm{loc}$ and $c_f$ can be explained by the boundary conditions: LFA is based on the assumption of periodic BC, and we use ABC.
Since the point source is located far away from the boundary, the actual convergence turns out to be slightly faster than the predictions, yet behaves similarly.
The LFA-predicted optimal damping parameters coincides with the numerically calculated ones. 
Therefore, in the remainder of this work we use these parameters for the fine grid relaxation.

\begin{figure}
\begin{center}
	\newcommand{\image}[1]{\includegraphics[width=0.3\linewidth]{#1}}
    \subfigure[\footnotesize Element patch]{\image{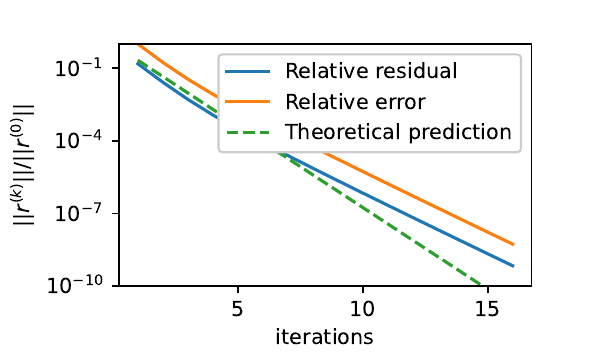}\label{fig:relres_4patch}} 
    \hspace{10pt}
    \subfigure[\footnotesize Plus patch]{\image{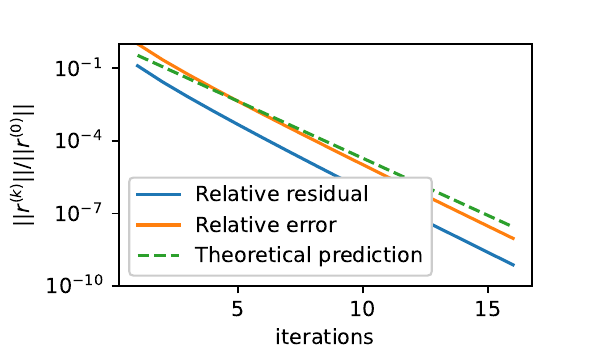}\label{fig:relres_5patch}}
    \hspace{10pt}
    \subfigure[\footnotesize RB patch]{\image{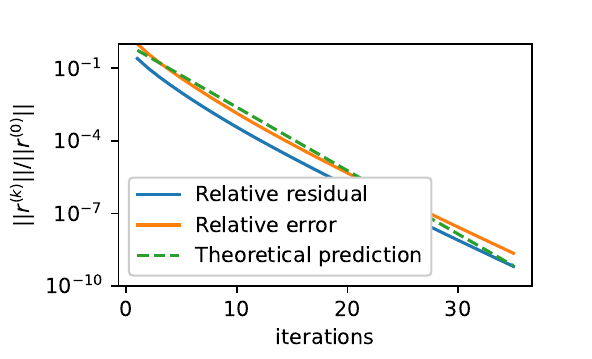}\label{fig:relres_skew5patch}}
\end{center}
\caption{Numerically calculated relative residual and relative error, and asymptotic expected convergence rate predicted by LFA for each 2D patch.
}\label{fig:relres}
\end{figure}

In Fig. \ref{fig:relres} we examine the asymptotic convergence behavior. 
In this experiment we solved a problem with a zero right-hand side (consequently, the exact solution is zero) and a random initial guess. 
In dashed straight lines we see the expected asymptotic behavior by LFA, and in solid curved line, the relative error ${\|\bfe^{(k)}\|}/{\|\bfe^{(0)}\|}$ and relative residual ${\|\bfr^{(k)}\|}/{\|\bfr^{(0)}\|}$.
The LFA convergence factors were taken as the minima from Fig. \ref{fig:damping}, and the dashed lines in Fig. \ref{fig:relres} represent a constant convergence rate equal to the asymptotic predicted rate.
The results show excellent agreement of the LFA predictions with the behavior in practice.

Judging by the two-grid predictions and numerical results, it seems that the RB patch is the worst additive Vanka patch. 
Surprisingly, increasing the number of level changes affects the convergence behavior, as shown in the results below.

\subsection{Designing the smoother and the coarse-grid correction}

First, we choose damping parameters for each smoother. 
In Fig. \ref{fig:damping} the damping parameters for the first level's relaxation were optimized separately by LFA for each of the smoothers (a similar experiment was conducted for damped Jacobi), and we observed by numerical experiments in Fig. \ref{fig:damping} that they are indeed optimal.
In the context of the Helmholtz equation, optimal damping parameters depend on the level of relaxation, since the corresponding number of grid points per wavelength differs for each coarse grid operator.
We optimize the damping parameters numerically as following: 
for the $l$-th level, we use exhaustive search on the damping parameter, assuming fixed damping parameters for the $1$-st to $(l-1)$-th level.
The resulting parameters for 2D and 3D are presented in Table \ref{tab:relaxParam}.
We observe that the damping parameters decrease with the level until the 3rd level relaxation and then increase.
This behavior is similar to the results in \cite{calandra2013improved}, where 2-level cycles were applied to systems with different number of grid points per wavelength with Jacobi relaxation.
The parameters in Table \ref{tab:relaxParam} suffice for up to 5-level methods.
In some of the experiments we used 6-level or 7-level methods, in which case we observed that the system is not very sensitive to the choice of the additional two parameters which are therefore not included in the table.

\begin{table}
\centering
\begin{tabular}{l|cccc|ccc}
\hline
  \toprule
  \mc{8}{c}{Damping parameters}\\
  \midrule
& \mc{4}{c|}{2D} & \mc{3}{c}{3D} \\
 & \small Jacobi & \small Element & \small $\,$Plus$\,$ & \small $\,\,\,\,$RB$\,\,\,\,$ & \small Jacobi & \small Element & \small $\,$Plus$\,$ \\
  \midrule
Fine level & 0.89 & 0.97 & 0.87 & 0.83 & 0.6 & 1.1 & 0.92 \\
First coarse level & 0.9 & 0.66 & 0.57 & 0.5 & 0.4 & 0.7 & 0.55 \\
Second coarse level & 0.3 & 0.48 & 0.55 & 0.4 & 0.3 & 0.45 & 0.45 \\
Third coarse level & 0.71 &  0.88 & 0.74 & 0.65 & 0.5 & 0.6 & 0.55 \\
  \bottomrule
 \end{tabular}
\caption{The damping parameters chosen by trial and error for different smoothers.
}
\label{tab:relaxParam}
\end{table}

\begin{table}
\centering
\begin{tabular}{c|ccc|ccc}
\hline
  \toprule
  \mc{7}{c}{GMRES iteration count for constant model}\\
  \midrule
 & \mc{3}{c|}{Jacobi, $\alpha=0.3$} & \mc{3}{c}{Element patch, $\alpha = 0.25$}  \\
\footnotesize Grid size (cells), freq.  & \footnotesize bicubic & \footnotesize mixed & \footnotesize lev-dep & \footnotesize bicubic & \footnotesize mixed & \footnotesize lev-dep \\
  \midrule
\footnotesize $128\times 128$, $\omega=25.6\pi$ & \footnotesize 39 (35) & \footnotesize 32 (30) & \footnotesize 35 (32) & \footnotesize 22 (21) & \footnotesize 23 (22) & \footnotesize 23 (23) \\
\footnotesize $256\times 256$, $\omega=51.2\pi$ & \footnotesize 70 (62) & \footnotesize 63 (55) & \footnotesize 64 (57) &  \footnotesize 40 (37) & \footnotesize 41 (38) & \footnotesize 42 (39) \\
\footnotesize $512\times 512$, $\omega=102.4\pi$ & \footnotesize 124 (111) & \footnotesize 111 (102) & \footnotesize 114 (102) & \footnotesize 74 (67) & \footnotesize 75 (69) & \footnotesize 75 (69) \\
  \midrule
  \midrule
 & \mc{3}{c|}{Plus patch, $\alpha = 0.25$} & \mc{3}{c}{RB patch, $\alpha = 0.18$} \\
 & \footnotesize bicubic & \footnotesize mixed & \footnotesize lev-dep & \footnotesize bicubic & \footnotesize mixed & \footnotesize lev-dep \\
  \midrule
\footnotesize $128\times 128$, $\omega=25.6\pi$ & \footnotesize 26 (25) & \footnotesize 26 (24) & \footnotesize 27 (25) & \footnotesize 20 (20) & \footnotesize 26 (25) & \footnotesize 20 (20) \\
\footnotesize $256\times 256$, $\omega=51.2\pi$ & \footnotesize 45 (41) & \footnotesize 55 (45) & \footnotesize 47 (43) & \footnotesize 36 (33) & \footnotesize 48 (44) & \footnotesize 36 (31) \\
\footnotesize $512\times 512$, $\omega=102.4\pi$ & \footnotesize 79 (71) & \footnotesize 100 (88) & \footnotesize 82 (74) & \footnotesize 64 (58) & \footnotesize 90 (85) & \footnotesize 63 (55) \\
  \bottomrule
 \end{tabular}
\caption{
GMRES(5) (in parentheses: non-restarted GMRES) iteration count for the solution of the 2D acoustic Helmholtz equation in homogeneous media, with a 4-level $W(1,1)$ cycle with different smoothers as a preconditioner. 
All the angular frequencies $\omega$ correspond to 10 grid points per wavelength.
}
\label{tab:const2D}
\end{table}

Second, we compare between different intergrid operators and different smoothers. Note that the two choices are not independent, as the coarse grid operator is defined by the Galerkin projection, so the choice of intergrid operators affects the computational stencil on which the smoother acts in coarser levels.
We compare the following intergrid schemes described in Subsection \ref{subsec:intergrid}:
\begin{itemize}
\item \textbf{Bicubic intergrid:} both $R$ and $P$ are calculated by \eqref{eq:bicubicIntergrid},
\item \textbf{Mixed intergrid:} $R$ is calculated by \eqref{eq:biliniearIntergrid} and $P$ is calculated by \eqref{eq:bicubicIntergrid}, and
\item \textbf{Level-dependent intergrid:} see \eqref{eq:level1Intergrid} and \eqref{eq:levelDepIntergrid}.
\end{itemize}
We solve the 2D Helmholtz equation in homogeneous media by GMRES(5) and by non-restarted GMRES, and use a $W(1,1)$ cycle as a preconditioner.
The results shown in Table \ref{tab:const2D} demonstrate that the level-dependent intergrid retains the efficiency of the bicubic intergrid operators, or even outperforms it, while presenting better sparsity and operator complexity.
Comparing the Element patch, Plus patch and RB patch for additive Vanka smoothing in Table \ref{tab:const2D}, we observe that the RB patch performs better for 4-level cycles, especially when combined with the level-dependent intergrid.
The shifts for each smoother are chosen by trial and error as the minimal shift that promises multigrid convergence\footnote{GMRES might converge even if the preconditioner is a divergent multigrid cycle. 
However, we observe that better results are achieved when the multigrid cycle converges.}.
This is done by calculating the convergence factor of the multigrid cycle for a relatively small grid, and from our experience a shift that promises convergence remains efficient for larger grids.
We observe that the results for non-restarted GMRES are lower by a few iterations.
However, we use GMRES(5) as our method of choice, for memory considerations.
Importantly, the RB patch requires a lower shift for convergence, and hence is expected to achieve better scalability.
Therefore, in the remainder of the 2D experiments we use the combination of level-dependent intergrid and RB patch.

This combination enables deep V-cycles with better scalability compared to existing methods, as shown below in Subsection \ref{subsec:scalability}.
We note that using deep W-cycles is overly expensive, due to the exponential growth in the number of coarse grid solutions, and for this reason our aim is to develop a method that will enable deep V-cycles. 
Nevertheless, we used ``shallow'' 4-levels method to choose the smoother for deeper cycles: 
shallow W-cycles give a better prediction for optimal smoother compared to shallow V-cycles, since in the W-cycle the coarse grids are visited more times.

\subsection{Scalability of the shift with respect to the level}
\label{subsec:scalability}

In this subsection we show a scalability property that our preconditioner achieves: scalability of the shift with respect to the level. 
In Fig. \ref{fig:shifts_256}, the shifts that optimize convergence are given.
We observe that using V-cycles with additive Vanka relaxation, a shift of no more than $\alpha=0.15$ suffices for any number of levels.
However, for V-cycles with Jacobi relaxation, the shift continues to increase, and reaches a plateau only at 6 levels, with a large shift of $\alpha=2.5$. 
We recall that the larger the shift, the method is less scalable with respect to the grid size.
To make the comparison of Vanka and Jacobi equitable in terms of computational cost, we included results for different number of pre- and post-relaxations.
And indeed, in Fig. \ref{fig:iterations_256}
 we see that deep V-cycles are more efficient with additive Vanka smoother than with damped Jacobi smoother, and this phenomenon does not depend on the grid resolution or the corresponding frequency.
By taking V-cycles of any depth, we enable a small and economical coarse grid problem, that needs to be solved only once.

\begin{figure}
\begin{center}
	\newcommand{\image}[1]{\includegraphics[width=0.4\linewidth]{#1}}
    \subfigure[\footnotesize Shifts]{\image{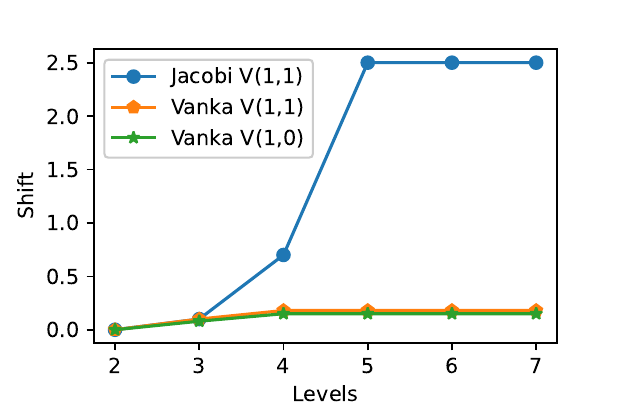}\label{fig:shifts_256}} 
    \hspace{20pt}
    \subfigure[\footnotesize Iterations]{\image{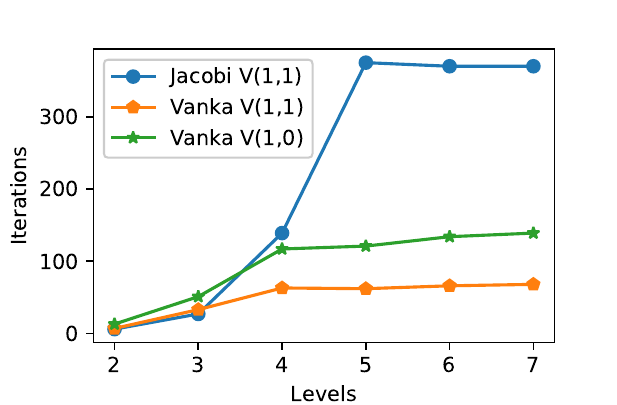}\label{fig:iterations_256}}
\\
\end{center}
\caption{Required shifts and GMRES(5) iteration count for the convergence of the 2D acoustic Helmholtz in homogeneous media, preconditioned by MG with Jacobi or RB patch additive Vanka as a smoother, for a grid of $256\times 256$ cells.
}
\label{fig:shift_scalability}
\end{figure}

\begin{remark}
Each Jacobi relaxation on the fine grid costs 10 $n_\mathrm{cells}$ FLOPs, 9 of them for the residual calculation.
Similarly, RB patch Vanka relaxation costs 34 $n_\mathrm{cells}$ FLOPs, about the same cost of 3 Jacobi relaxations.
In each of the coarse levels, however, the residual calculation costs 25 $n_\mathrm{cells}$ FLOPs, and hence RB patch Vanka costs only twice as Jacobi costs. 
It would be reasonable then to compare Vanka V(1,0) with Jacobi V(1,2).
Nevertheless, we observe that Jacobi V(1,1) in fact performs better than Jacobi V(1,2), since Jacobi is not stable for indefinite problems, and additional relaxations might lead to divergence.
\end{remark}

\subsection{Performance for heterogeneous media}

\begin{figure}
\begin{center}
	\newcommand{\image}[1]{\includegraphics[width=0.2\linewidth]{#1}}
    \subfigure[\footnotesize $\kappa^2(\bfx)$, linear]{\image{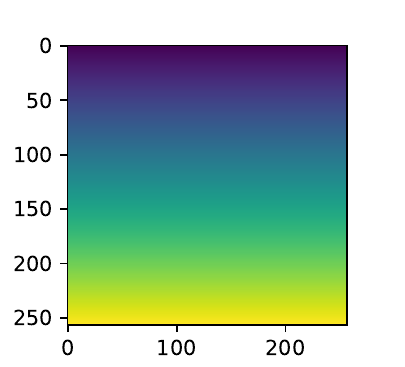}\label{fig:Linear_square}} 
    \hspace{5pt}
    \subfigure[\footnotesize $p(\bfx)$, linear]{\image{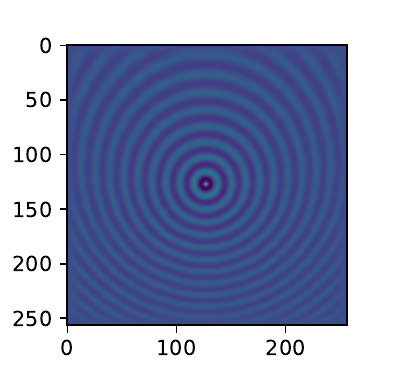}\label{fig:waveform_linear}} 
    \hspace{5pt}    
    \subfigure[\footnotesize $\kappa^2(\bfx)$, wedge]{\image{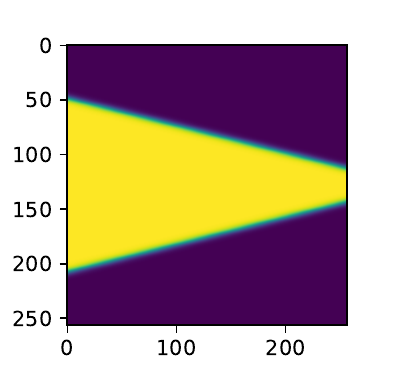}\label{fig:Wedge_square}}
        \hspace{5pt}
    \subfigure[\footnotesize $p(\bfx)$, wedge]{\image{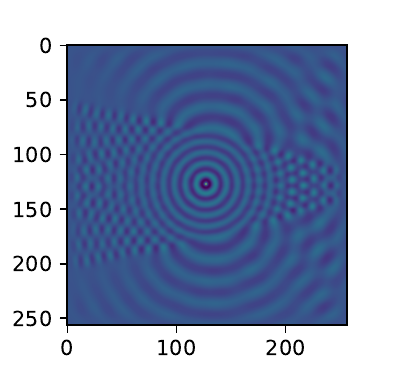}\label{fig:waveform_wedge}} 
\\
\end{center}
\caption{The slowness model and resulting waveform for linear and wedge heterogeneous velocity models in 2D.
}\label{fig:models}
\end{figure}

In this section we demonstrate the performance of our method for the following heterogeneous velocity models:
\begin{itemize}
\item \textbf{Linear:} a dimensionless domain $[0,1]^2$, on which the slowness squared varies linearly in the vertical direction in the range $\kappa^2(x)\in [0.25,1]$, see Fig. \ref{fig:Linear_square}\, or the range $\kappa^2(x)\in [0.05,1]$.
\item \textbf{Wedge:} a dimensionless domain $[0,1]^2$ on which $\kappa^2(x)\in [0.25,1]$ varied with a rapid change between two constant velocity domains, see Fig. \ref{fig:Wedge_square}.
We also use the same model with the range $\kappa^2(x)\in [0.05,1]$.
\end{itemize}
Both models have the same variation ranges; one varies smoothly, and the other rapidly.

\begin{table}
\centering
\begin{tabular}{c|ccc|ccc}
\hline
  \toprule
  \mc{7}{c}{GMRES(5) iteration count for 2D heterogeneous media}\\
  \midrule
 & \mc{3}{c|}{Linear, $\kappa^2\in[0.25,1]$} & \mc{3}{c}{Wedge, $\kappa^2\in[0.25,1]$}  \\
& 2-level & 3-level & 4-level & 2-level & 3-level & 4-level \\
Grid size (cells), freq. & \footnotesize $\alpha=0.0$ & \footnotesize $\alpha=0.1$ & \footnotesize $\alpha=0.25$ & \footnotesize $\alpha=0.0$ & \footnotesize $\alpha=0.1$ & \footnotesize $\alpha=0.25$\\
  \midrule
\footnotesize $128\times 128$, $\omega=25.6\pi$ & 6 & 12 & 34 & 6 & 19 & 42 \\
\footnotesize $256\times 256$, $\omega=51.2\pi$ & 6 & 20 & 76 & 6 & 31 & 71 \\
\footnotesize $512\times 512$, $\omega=102.4\pi$ & 6 & 35 & 182 & 7 & 53 & 137 \\
\footnotesize $1024\times 1024$, $\omega=204.8\pi$ & 6 & 63 & 512 & 7 & 100 & 275 \\
\midrule
\midrule
 & \mc{3}{c|}{Linear, $\kappa^2\in[0.05,1]$} & \mc{3}{c}{Wedge, $\kappa^2\in[0.05,1]$} \\
& 2-level & 3-level & 4-level & 2-level & 3-level & 4-level \\
Grid size (cells), freq. & \footnotesize $\alpha=0.0$ & \footnotesize $\alpha=0.1$ & \footnotesize $\alpha=0.25$ & \footnotesize $\alpha=0.0$ & \footnotesize $\alpha=0.1$ & \footnotesize $\alpha=0.3$\\
  \midrule
\footnotesize $128\times 128$, $\omega=25.6\pi$ & 6  & 12 & 32 & 6 & 22 & 65 \\
\footnotesize $256\times 256$, $\omega=51.2\pi$ & 6 & 21 & 70 & 7 & 36 & 116 \\
\footnotesize $512\times 512$, $\omega=102.4\pi$ & 6 & 37 & 161 & 7 & 66 & 224 \\
\footnotesize $1024\times 1024$, $\omega=204.8\pi$ & 6 & 67 & 445 & 7 & 123 & 436 \\
  \bottomrule
 \end{tabular}
\caption{GMRES(5) iteration count for the solution of the 2D acoustic Helmholtz equation in different heterogeneous media.
A $V(1,1)$ cycle with RB patch additive Vanka smoother and level-dependent intergrid operators serves as a preconditioner.
All the angular frequencies $\omega$ correspond to 10 grid points per wavelength.
}
\label{tab:heterogeneous2D}
\end{table}

We solve the 2D Helmholtz equation in linear and wedge model with different heterogeneity ranges using GMRES(5) preconditioned by multigrid cycles with level-dependent intergrid and RB patch additive Vanka smoother.
We use the damping parameters reported in Table \ref{tab:relaxParam}.
The results in Table \ref{tab:heterogeneous2D} show similar scalable behavior for 2-levels cycles for both velocity models.
There are some differences in the convergence rate for more than 2 levels: for 3-level cycles, the linear model of range $\kappa^2\in[0.25,1]$ converge faster than the wedge model of the same range, and for the 4-level cycles, the wedge model of range $\kappa^2\in[0.25,1]$ requires less shift than the linear model of the same range.
For the larger range, $\kappa^2\in[0.05,1]$, the linear and wedge models behave similarly in 4 levels.
Although there are some differences in iteration count that correlates naturally to the heterogeneity range and contrast, the differences between the models are not decisive, and compared to the rightmost column of Table \ref{tab:const2D}, it is evident that the efficiency of the method is not strongly affected by heterogeneity.

In Fig. \ref{fig:shift_sensitivity} we depict the iteration count vs. the shift for the linear and wedge model of different ranges.
We observe that the sensitivity of the Vanka-smoothed cycle to the complex shift is lower compared to the sensitivity of damped Jacobi.

\begin{figure}
\begin{center}
	\newcommand{\image}[1]{\includegraphics[width=0.48\linewidth]{#1}}
    \subfigure[\footnotesize RB Vanka]{\image{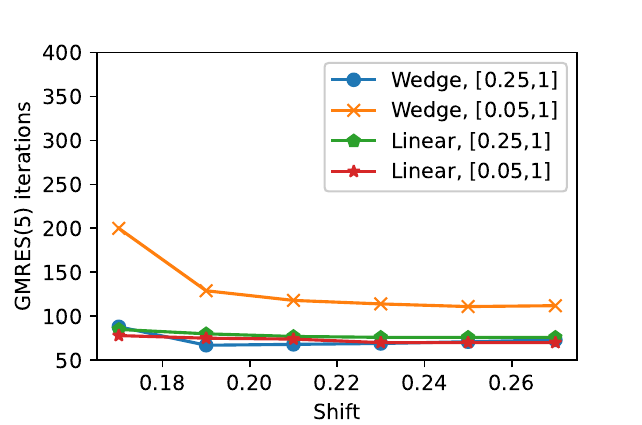}\label{fig:shift_sensitivity_Vanka}}
    \hspace{5pt}
    \subfigure[\footnotesize Jacobi]{\image{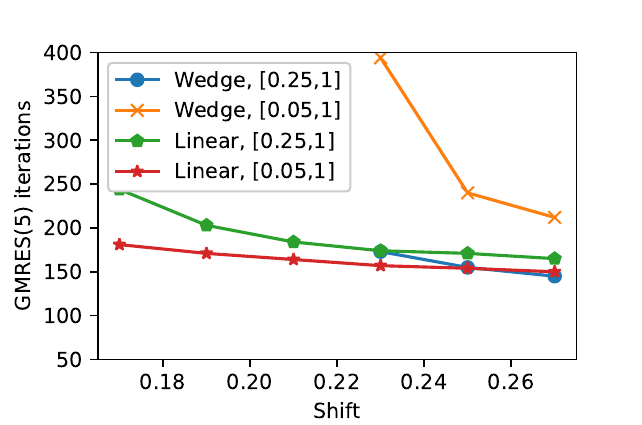}\label{fig:shift_sensitivity_Jacobi}}\\
\end{center}
\caption{GMRES(5) iteration count for solving the Helmholtz equation in linear or wedge media with different ranges, for grid size $256\times256$.
A multigrid $V(1,1)$ cycle of 4 levels with RB patch additive Vanka or damped Jacobi serves as a preconditioner.
For both ranges of the wedge model, three points are missing from the graph on the right, since the iteration count exceeded 1000.
}\label{fig:shift_sensitivity}
\end{figure}

\subsection{3D experiments}

In this subsection we demonstrate the performance of our method in 3D experiments. 
To this end, we first choose appropriate components for the multigrid cycle. 
Some of them are inherited from the developed 2D method, e.g., the computational advantage of the level-dependent scheme over tricubic intergrid operators is evident and even more significant in the 3D case: 
the resulting Galerkin coarse stencils are $5\times5\times5$ instead of $7\times7\times7$.
Nonetheless, its gain might not pay off compared to the trilinear intergrid operators, which have stencils of only $3\times3\times3$ points.

However, designing the additive Vanka smoother is more cumbersome:
clearly, the 27 points Full patch from Fig. \ref{fig:27patch} is too expensive to consider. 
Even the 3D RB patch from Fig. \ref{fig:13patch}, the analogue of the preferred 2D patch, would be considerably more expensive than the Element patch from Fig. \ref{fig:8patch} or Plus patch from Fig. \ref{fig:7patch}, since it has 13 points compared to only 7-8 points.
Not only that the RB patch will have to perform at least twice better to justify its size, it also loses the advantage of quick inversion of very small patch-matrices, possible in some implementations.
Hence, to keep a reasonable computational cost per iteration, we only compare the Element patch and the Plus patch.
The damping parameters used for each smoother at each level are found by trial and error, and are reported in Table \ref{tab:relaxParam}. 
The same damping parameters are used for homogeneous and heterogeneous media.

The results of the comparison between these patches are presented in Table \ref{tab:smoother-choice3D}. 
For reference, we include results for damped Jacobi as well.
We compare the performance of the three smoothers for relatively small grid sizes of homogeneous media. 
We use 5-level $V(1,1)$-cycles with level-dependent intergrid.
The shifts are chosen by trial and error to be the minimal shifts that promise convergence of the multigrid for all of the examined grid sizes. 
The Jacobi relaxation leads to significant increase of iteration counts with the grid size and frequency, as well as the Plus patch Vanka smoother.
The Element patch leads to smaller iteration count and significantly less computation time.
This result fits the physical behavior of solutions of the 3D Helmholtz equation: 
the sphere-shaped pattern around the source is better captured by a more compact patch, such as the Element patch.

\begin{table}
\centering
\begin{tabular}{l|lll}
\hline
  \toprule
  \mc{4}{c}{GMRES(5) iteration count and time (sec) for 3D constant media}\\
  \midrule
 & Jacobi & Element & Plus \\
Grid size (cells), freq. & \footnotesize $\alpha = 0.5$ & \footnotesize $\alpha = 0.4$ &  \footnotesize $\alpha = 0.65$ \\
  \midrule 
   $48\times 48 \times 48$, $\omega=9.6\pi$ & 21 (1.36) & 13 (1.29) & 18 (1.75) \\
 $64\times 64 \times 64$, $\omega=12.8\pi$ & 35 (8.30) & 20 (5.05) & 22 (5.40) \\
 $96\times 96 \times 96$, $\omega=19.2\pi$ & 90 (53.74) & 34 (62.53) & 72 (136.52) \\
 $128\times 128 \times 128$, $\omega=25.6\pi$ & 143 (371.47) & 49 (198.00) & 138 (656.37) \\
  \bottomrule
 \end{tabular}
\caption{
GMRES(5) iteration count for the solution of the 3D acoustic Helmholtz equation in homogeneous media, with a 5-level $V(1,1)$ cycle with different smoothers as a preconditioner, and level-dependent intergrid operators, for different grid sizes and corresponding angular frequencies for 10 grid points per wavelength. 
In parentheses: wall-clock time in seconds.
}
\label{tab:smoother-choice3D}
\end{table}

In the second experiment, presented in Fig. \ref{fig:time_all}, we hold time comparison.
Each reported value is the average of 6 runs of the same experiment.
First, in Fig. \ref{fig:time} we compare our method (Element patch additive Vanka smoother and level-dependent intergrid) to the ``plain'' shifted Laplacian multigrid method (with damped Jacobi smoother and trilinear intergrid), using 5-level $V(1,1)$-cycles for both methods.
We observe that the our method performs better: for the largest grid examined, the convergence is approximately $2.18$ times faster.
We emphasize that the improvement is not only quantitative, but also qualitative: the scaling is improved, leading to nearly linear time growth.
The better scaling can be measured by a smaller growth in average: 
about $78$ seconds per million unknowns in our method, compared to about $162$ seconds of time growth using damped Jacobi trilinear intergrid.
Second, we aim to explore the scalability with respect to the level, that is, to see if the properties from Fig. \ref{fig:shift_scalability} are preserved in 3D.
To this end, we compare Vanka to Jacobi, both with level-dependent intergrid, for different number of levels and a constant grid of $128\times128\times128$ cells. 
The results in Fig. \ref{fig:time_levels} show that scalability with respect to the number of levels is retained. 
The graph in Fig. \ref{fig:time_levels} does not include results for a 2-level method, which turned out to be almost an order of magnitude slower because of the fill-in of the large LU factors.
Results for a 7-level method include only Vanka since Jacobi failed to converge after 1000 iterations.
We note that the shift choice for a large number of levels in 3D was always $\alpha=0.5$.
For Vanka, this shift was chosen since it was the smallest that led to a convergent multigrid cycle (similarly to the shift choice strategy in 2D), and hence the improved scalability holds. 
For Jacobi, the minimal shifts that led to convergent cycles in 5--7 level methods were so large, so that the preconditioner was no longer efficient. 
The shift $\alpha=0.5$ was found to be the best compromise in practice.

\begin{figure}
\begin{center}
	\newcommand{\image}[1]{\includegraphics[width=0.48\linewidth]{#1}}
    \subfigure[\footnotesize Time vs. grid size]{\image{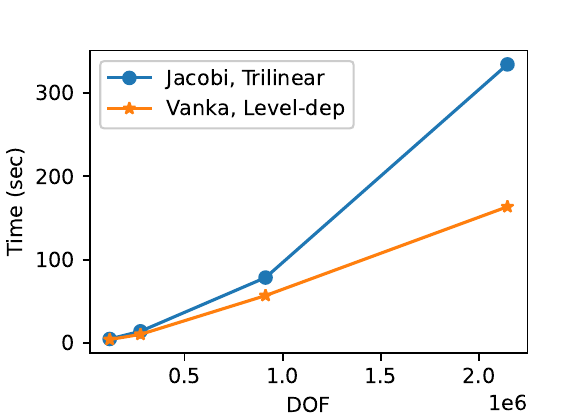}\label{fig:time}}
    \hspace{5pt}
    \subfigure[\footnotesize Time vs. levels]{\image{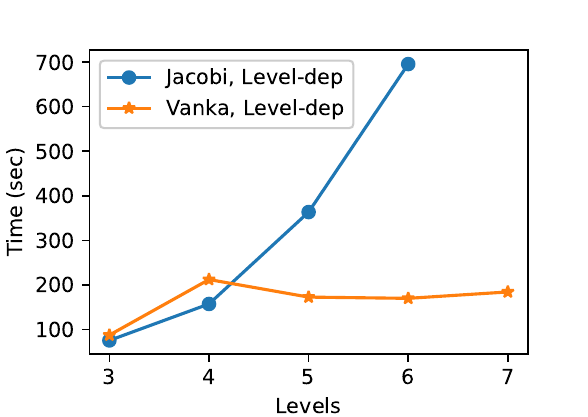}\label{fig:time_levels}}\\
\end{center}
\caption{Average wall-clock time (sec) for the solution of the 3D acoustic Helmholtz equation in linear media with $\kappa^2\in[1,2]$ using GMRES(5), preconditioned by $V(1,1)$ cycles.
On the left, 5-level cycles for different grid sizes, and on the right, for a $128\times128\times128$ cells grid and different number of levels. 
The shift is $\alpha=0.2$ for a 3-level method and $\alpha=0.5$ for larger number of levels.
}\label{fig:time_all}
\end{figure}

\begin{remark}\label{rem:setup}
The time comparisons in Fig. \ref{fig:time_all} compare only the running time of the GMRES(5) iteration, and not the setup time of the algorithm.
Clearly, Vanka's setup includes more FLOPs, as it includes inverting many small patch-matrices. 
Naturally, it might lead to larger setup times.
However, the setup time depends heavily on the implementation and on how the sparse matrices are stored.
Moreover, in many applications --- e.g., full waveform inversion and scattering with multiple sources --- the same Helmholtz operator is solved against hundreds or even thousands of different right hand sides.
Therefore, the differences in setup time become negligible in comparison to the overall cost.
\end{remark}

In the third experiment we demonstrate the performance of our method for the Overthrust model \cite{aminzadeh19973}, a real-world geophysical velocity model with jumping coefficients, see Fig. \ref{fig:Overthrust}.
The model is shallow, so we extend it vertically by 12 to 16 grid points, to facilitate enforcing the absorbing boundary layer. 
In Fig. \ref{fig:time_levels_overthrust} we compare the average wall-clock time (averaged over 6 experiments) needed for the solution of the 3D Helmholtz equation by GMRES(5), preconditioned by a multigrid V(1,1) cycle with trilinear intergrid operators.
We chose standard trilinear intergrid operators for memory reasons, and also since we observe that in large 3D cases, the gain of the wide intergrid operators does not justify its high computational cost.
We compare three smoothers: damped Jacobi, Element patch additive Vanka, and a combination of Jacobi on the fine level and Vanka on the rest of the levels. 
In Fig. \ref{fig:time_levels_overthrust} we hold this comparison for a grid of $256\times256\times96$ cells, and for 3, 4, 5 and 6 levels.
We observe that for 5 levels and more, Vanka outperforms Jacobi. 
Nevertheless, the combined method outperforms both in all the examined numbers of levels, and scales better with respect to the levels.

\begin{figure}
  \centering
  \includegraphics[width=0.6\textwidth]{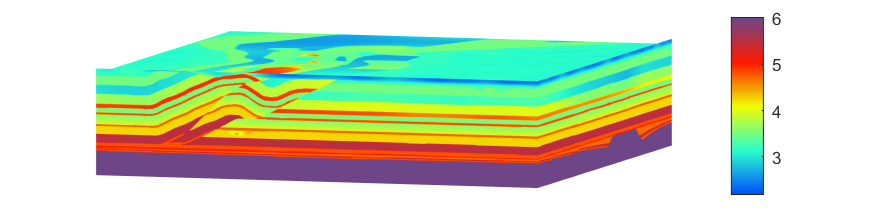}\\
  \caption{The Overthrust model. Velocity units: $km/sec$.}\label{fig:Overthrust}
\end{figure}

\begin{figure}
  \centering
  \includegraphics[width=0.5\textwidth]{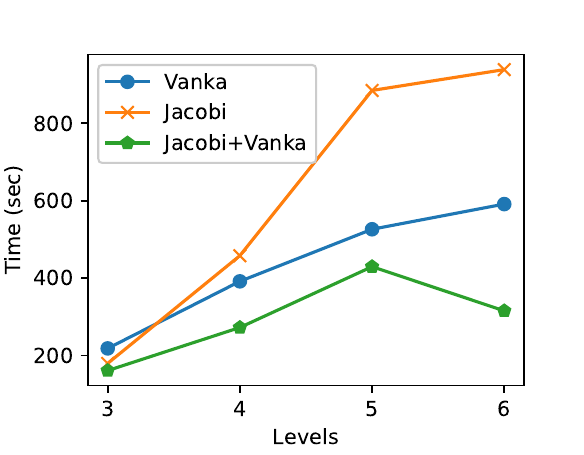}\\
  \caption{Wall-clock time needed for solving the Helmholtz equation in Overthrust media for a grid of $256\times256\times96$ cells, for different number of levels, using preconditioned GMRES(5).
A multigrid $V(1,1)$-cycle serves as a preconditioner. 
We used three options of smoothers: Element patch additive Vanka, damped Jacobi relaxation, or Jacobi on the fine level and Vanka on the coarse level, and trilinear intergrid operators.}\label{fig:time_levels_overthrust}
\end{figure}

\begin{table}
\centering
\begin{tabular}{c|ccc}
\hline
  \toprule
  \mc{4}{c}{GMRES(5) iteration count and time (sec) for Overthrust media}\\
  \midrule
 &  \mc{3}{c}{Jacobi+Vanka} \\
 & \footnotesize 3-level &  \footnotesize 4-level & \footnotesize 5-level  \\
Grid size (cells), freq. & \footnotesize $\alpha = 0.2$ & \footnotesize $\alpha = 0.4$ & \footnotesize $\alpha = 0.5$ \\
  \midrule
\footnotesize $128\times 128\times 56$, $\omega = 2.86 \pi$ & \footnotesize 12 (10.60) & \footnotesize 19 (15.59) & \footnotesize  NA \\
\footnotesize $192\times 192\times 72$,  $\omega = 4.22 \pi$ & \footnotesize 15 (43.90) & \footnotesize 30 (67.77) & \footnotesize  NA \\
\footnotesize $256\times 256\times 96$, $\omega = 5.58 \pi$ & \footnotesize 20 (160.81) & \footnotesize 45 (272.16) & \footnotesize 69 (429.39) \\
\footnotesize $320\times 320\times 112$, $\omega = 6.97 \pi$ & \footnotesize 28 (409.03) & \footnotesize 79 (735.80) & \footnotesize 120 (1084.76) \\
\midrule
\midrule
 & \mc{3}{c}{Jacobi} \\
  \midrule
\footnotesize $128\times 128\times 56$, $\omega = 2.86 \pi$ & \footnotesize 13 (13.65) & \footnotesize 25 (22.96) & \footnotesize  NA \\
\footnotesize $192\times 192\times 72$,  $\omega = 4.22 \pi$ & \footnotesize 16 (53.04) & \footnotesize 41 (114.51) & \footnotesize  NA \\
\footnotesize $256\times 256\times 96$, $\omega = 5.58 \pi$ & \footnotesize 22 (179.30) & \footnotesize 76 (457.86) & \footnotesize 173 (884.90) \\
\footnotesize $320\times 320\times 112$, $\omega = 6.97 \pi$ & \footnotesize 30 (425.34) & \footnotesize 174 (1635.95) & \footnotesize  $>1000$ \\
  \bottomrule
 \end{tabular}
\caption{GMRES(5) iteration count for the solution of the 3D acoustic Helmholtz equation in Overthrust media for different grid sizes and the corresponding angular frequencies with 10 grid points per wavelength.
A $V(1,1)$ cycle with damped Jacobi relaxation on the fine level and Element patch additive Vanka smoother on the rest of the levels, compared to damped Jacobi smoother on all of the levels, both with trilinear intergrid operators, serves as a preconditioner.
In parentheses: wall-clock time in seconds.
NA stands for ``not available'', when the number of cells does not enable a 5-level method.
}
\label{tab:Overthrust}
\end{table}

In Table \ref{tab:Overthrust}
we measure the GMRES(5) iterations and the average time (averaged over 6 runs) needed for the solution of the 3D acoustic Helmholtz problem in Overthrust media, using 3-level and 4-level multigrid methods with the combination of damped Jacobi on the fine level and Element patch additive Vanka on the rest of the levels as a preconditioner, compared to cycles with damped Jacobi.
For both methods, we used standard trilinear intergrid operators, as before, for memory and computational cost considerations.
The results in Table \ref{tab:Overthrust} establish the efficiency of the combined smoothing method for real 3D cases, and suggest that the Vanka smoother can enhance performance significantly for highly heterogeneous velocity models.
Interestingly, both in Table \ref{tab:Overthrust} and in Fig. \ref{fig:time_levels}, the 3-level method is the fastest. 
It can be explained by the lower shift that a 3-level method requires.
Nevertheless, the advantage of deeper V-cycles becomes apparent in extremely large cases, where even the coarse grid of the 3-level method is too large to solve exactly.
To summarize, we suggest a 3-level method with combined Jacobi-Vanka smoother whenever the method is possible, and when memory limitations makes a 3-level method infeasible, one can increase the levels without a large degradation in the computational cost.

\section{Conclusion}\label{sec:conclusion}

In this work, we introduced a Vanka-smoothed shifted Laplacian preconditioner for the acoustic Helmholtz equation. By a carefully designed additive Vanka smoother we achieve better scalability. 
For 2D, we also equip the method with a level-dependent intergrid scheme that enables a lower complex shift while keeping the operator complexity reasonable. 
We utilize two-grid LFA to tune the damping parameters, and show that the theoretical predictions are well-aligned with the convergence in practice.
However, we observe that the two-grid behavior is not a good predictor of the preferred smoother for deeper cycles, which enables a small coarsest grid problem and saves time and memory.
For cycles of at least 4 levels, we observe that the RB patch additive Vanka smoother performs best in 2D, and the Element patch additive Vanka smoother performs best in 3D.
For large and heterogeneous 3D problems, a combined approach of damped Jacobi on the fine level and Element patch additive Vanka on the rest of the levels was found to be more beneficial than Vanka smoothing on all of the levels.

We observed that 3-level cycles has the fastest convergence for the examined problems, probably since they require a lower complex shift. 
Nonetheless, since the complexity of the LU factors of the coarse grid grow super-linearly, for extremely large problems we expect that 3-level methods will not be applicable. 
Hence, using more levels without significantly degrading convergence is a sought property. 
Importantly, 
 deep V-cycles are convergent using our method. 
It enables a faster and more robust solution, and in practice, we observe a nearly linear time scaling in 3D, and an improved scaling even for challenging geophysical media.

\section*{Acknowledgment} The authors thank Chen Greif for constructive comments and insightful discussions that provided valuable guidance during the preparation of this work and for helping us clarify and articulate our ideas.

\bibliographystyle{siamplain}

\bibliography{main.bbl}

\end{document}